\documentclass[twoside,reqno]{amsart}

\usepackage{latexsym}

\newcommand{\qdn}{\hspace*{-1.5mm}}
\newcommand{\qqdn}{\hspace*{-2.5mm}}
\newcommand{\xqdn}{\hspace*{-5.0mm}}
\newcommand{\xxqdn}{\hspace*{-10mm}}

\newcommand{\fns}{\footnotesize}
\newcommand{\sst}{\scriptstyle}
\newcommand{\sss}{\scriptscriptstyle}

\newcommand{\fnk}[3]{\left[\qdn\ba{#1}#2\\#3\ea\qdn\right]}
\newcommand{\ffnk}[4]{\left[\qdn\ba{#1}#3\\#4\ea{\!\Big|\:#2}\right]}

\newcommand{\binm}{\binom}

\newcommand{\nnm}{\nonumber}
\newcommand{\be}{\begin{equation}}
\newcommand{\ee}{\end{equation}}
\newcommand{\ba}{\begin{array}}
\newcommand{\ea}{\end{array}}
\newcommand{\bmn}{\begin{eqnarray}}
\newcommand{\emn}{\end{eqnarray}}
\newcommand{\bnm}{\begin{eqnarray*}}
\newcommand{\enm}{\end{eqnarray*}}
\newcommand{\bln}{\begin{subequations}}
\newcommand{\eln}{\end{subequations}}
\newcommand{\blm}{\bln\bmn}
\newcommand{\elm}{\emn\eln}

\newtheorem{thm}{Theorem}
\newtheorem{lemm}[thm]{Lemma}
\newtheorem{corl}[thm]{Corollary}
\newtheorem{prop}[thm]{Proposition}

\newtheorem{entry}{Entry}

\newcommand{\bbtm}[4]{\bibitem{kn:#1}{#2,}~\emph{#3,}~{#4.}}
\newcommand{\cito}[1]{\cite{kn:#1}}
\newcommand{\citu}[2]{\cite[#2]{kn:#1}}

\begin{document} 
{\fns 
\title{Contiguous relations of $_3\phi_2$-series}
\author{$^a$Chuanan Wei and $^b$Dianxuan Gong}
 \dedicatory{ $^A$Department of Information Technology\\
            Hainan Medical College, Haikou 571199, China\\
           $^B$College of Sciences\\
 Hebei United University, Tangshan 063009, China}
\thanks{\emph{Email addresses}: weichuanan@yahoo.com.cn (C. Wei),
gongdianxuan@yahoo.com.cn (D. Gong)}
\address{ }
\footnote{\emph{2010 Mathematics Subject Classification}: Primary
33D15 and Secondary 11J72.}
 \keywords{Abel's lemma;
 Two-term contiguous relation of $_3\phi_2$-series;
 Three-term contiguous relation of $_3\phi_2$-series}

\begin{abstract}
According to Abel's lemma and the method of linear combinations, we
establish numerous contiguous relations of $_3\phi_2$-series, which
can be regarded as $q$-analogues of the contiguous relations of
$_3F_2$-series due to Krattenthaler and Rivoal \cito{krattenthaler}
or Chu and Wang \cito{wang}.
\end{abstract}

\maketitle\thispagestyle{empty}
\markboth{C. Wei and D. Gong}
         {Contiguous relations of $_3\phi_2$-series}


\section{Introduction}
For two complex numbers $x$ and $q$ with $|q|<1$, define the
$q$-shifted factorial by
 \bnm
 (x;q)_n=
\begin{cases}
\prod_{i=0}^{n-1}(1-xq^i),&\quad n>0;\\
1,&\quad n=0;\\
\frac{1}{\prod_{j=n}^{-1}(1-xq^j)},&\quad n<0.
\end{cases}
 \enm
The fractional form of it reads as
\[\ffnk{cccc}{q}{\alpha,&\beta,&\cdots,&\gamma}{A,&B,&\cdots,&C}_n=
\frac{(\alpha;q)_n(\beta;q)_n\cdots(\gamma;q)_n}
{(A;q)_n(B;q)_n\cdots(C;q)_n}.\]

 Following Gasper and
Rahman~\cito{gasper}, the basic hypergeometric series can be defined
by
\[\xqdn_{r}\phi_s\ffnk{cccc}{q;z}{a_1,a_2,\cdots,a_r}{b_1,b_2,\cdots,b_s}
 =\sum_{k=0}^\infty
\frac{(a_{1};q)_{k}(a_2;q)_{k}\cdots(a_r;q)_{k}}
 {(b_1;q)_{k}(b_2;q)_{k}\cdots(b_s;q)_{k}}\Big\{(-1)^kq^{\binm{k}{2}}\Big\}^{1+s-r}\frac{z^k}{(q;q)_k},\]
where $\{a_{i}\}$ and $\{b_{j}\}$ are complex parameters such that
no zero factors appear in the denominators of the summand on the
right hand side. Throughout the paper, we shall also use the
\emph{shifted}-basic hypergeometric series
\[_{r}\phi_s^{*}\ffnk{cccc}{q;z}{a_1,a_2,\cdots,a_r}{b_1,b_2,\cdots,b_s}
 =\sum_{k=0}^\infty
\frac{(a_{1};q)_{k}(a_2;q)_{k}\cdots(a_r;q)_{k}}
 {(b_1;q)_{k}(b_2;q)_{k}\cdots(b_s;q)_{k}}\Big\{(-1)^kq^{\binm{k}{2}}\Big\}^{1+s-r}\frac{(1\!-\!q^k)z^k}{(q;q)_k},\]
whose summation index begins essentially with $k=1$, instead of
$k=0$.

 For a complex sequence $\{\tau_k\}$,
define respectively the forward difference operator
$\tilde{\triangle}$ and the backward difference operator $\nabla$ by
\[\tilde{\triangle}\tau_k=\tau_k-\tau_{k+1} \quad\text{and}\quad \nabla\tau_k=\tau_k-\tau_{k-1}. \]
 Then Abel's lemma (cf. \cito{wang}) can be
stated as follows.

\begin{lemm}
\label{lemm}
 For two complex sequences $\{U_k\}$ and
$\{V_k\}$, there holds the relation:
\[\sum_{k=0}^{\infty}U_k{\tilde{\triangle} V_k}=\sum_{k=0}^{\infty}V_k{\nabla U_k}\]
provided that one of the series on both sides converges,
$U_{-1}V_0=0$ and $U_kV_{k+1}\to0$ as $k\to \infty$.
\end{lemm}

 There are many contiguous relations in the literature. Several interesting ones can
be seen in the papers \cito{buschman}-\cito{wang},
\cito{gupta-a}-\cito{krattenthaler}, \cito{purohit}-\cito{rakha-c}
and \cito{swarttouw}-\cito{vidunas-b}. Implied by the work just
mentioned, we shall give numerous contiguous relations of
$_3\phi_2$-series in terms of Lemma \ref{lemm} and the method of
linear combinations.

The structure of this paper is arranged as follows. In section 2, we
shall use the Lemma \ref{lemm} to found four three-term contiguous
relations of $_3\phi_2$-series denominated by patterns \textbf{A},
\textbf{B}, \textbf{C} and \textbf{D}. Then they will be applied to
offer numerous two- and three-term contiguous relations of
$_3\phi_2$-series in sections 3-4 in accordance with the method of
linear combinations.
\section{Four three-term contiguous relations of $_3\phi_2$-series}
In this section, we show four patterns \textbf{A}, \textbf{B},
\textbf{C} and \textbf{D} satisfied by three
 $_3\phi_2$-series through Abel's lemma. Throughout this section, we assume that the
parameters of all the $_3\phi_2$-series are subject to the condition
$|bd/qace|<1$ in order that Lemma \ref{lemm} can be applied
smoothly.
\subsection{Pattern A}$ $\\\\
Define two sequences by
\[U_k=\ffnk{ccc}{q}{qa,&d/a}{q,&d}_k \quad\text{and}\quad V_k=\ffnk{ccc}{q}{c,&e}{b,&d/qa}_k\bigg(\frac{bd}{qace}\bigg)^k.\]
Then it is not difficult to check the limiting relation
\[U_{-1}V_{0}=\lim_{n\to \infty}U_{n}V_{n+1}=0\]
and the finite differences
 \bnm
&&\xqdn\nabla U_k=\ffnk{ccc}{q}{a,&d/qa}{q,&d}_kq^k,\\
&&\xqdn\tilde{\triangle} V_k=\ffnk{ccc}{q}{c,&e}{qb,&d/a}_k
 \bigg(\frac{bd}{ace}\bigg)^k
 \bigg\{\frac{\sss(1-b)(1-d/qa)-(1-c)(1-e)\frac{bd}{qace}}{(1-b)(1-d/qa)}
 +\frac{(1-bd/qace)}{(1-b)(1-d/qa)}\frac{1-q^k}{q^k}\bigg\}.
 \enm
In accordance with Lemma \ref{lemm}, we can manipulate the following
$_3\phi_2$-series:
 \bnm
&&_3\phi_2\ffnk{cccc}{q;\frac{bd}{ace}}{a,c,e}{b,d}
=\sum_{k\geq0}V_k\nabla U_k=\sum_{k\geq0}U_k\tilde{\triangle} V_k\\
&&\:\:=\:\:\frac{\sss(1-b)(1-d/qa)-(1-c)(1-e)\frac{bd}{qace}}{(1-b)(1-d/qa)}
\sum_{k\geq0}\ffnk{ccc}{q}{qa,&c,&e}{q,&qb,&d}_k\bigg(\frac{bd}{ace}\bigg)^k\\
 &&\:\:+\:\:\:\frac{(1-bd/qace)}{(1-b)(1-d/qa)}
 \sum_{k\geq0}(1-q^k)\ffnk{ccc}{q}{qa,&c,&e}{q,&qb,&d}_k
 \bigg(\frac{bd}{qace}\bigg)^k.
 \enm
Shifting the summation index $k\to k+1$ for the last sum, we obtain
the following relation.

\begin{thm}[Pattern A]\label{thm-a} For five complex numbers $\{a, b,
c, d, e\}$ subject to the condition $|bd/qace|<1$, there holds the
three-term contiguous  relation of $_3\phi_2$-series:
 \blm
&&_3\phi_2\ffnk{cccc}{q;\frac{bd}{ace}}{a,c,e}{b,d}
=\label{eq-a}\mathcal{A}_q{_3\phi_2}\ffnk{cccc}{q;\frac{bd}{ace}}{qa,c,e}{qb,d}
+\mathbb{A}_q{_3\phi_2}\ffnk{cccc}{q;\frac{bd}{qace}}{q^2a,qc,qe}{q^2b,qd},\\
&&_3\phi_2\ffnk{cccc}{q;\frac{bd}{ace}}{a,c,e}{b,d}
=\label{eq-aa}\mathcal{A}_q{_3\phi_2}\ffnk{cccc}{q;\frac{bd}{ace}}{qa,c,e}{qb,d}
+\mathfrak{A}_q{_3\phi_2}^{*}\ffnk{cccc}{q;\frac{bd}{qace}}{qa,c,e}{qb,d},
 \elm
where the coefficients $\mathcal{A}_q$, $\mathbb{A}_q$ and
$\mathfrak{A}_q$ are defined by
 \bnm
&&\mathcal{A}_q:=\mathcal{A}_q(a,c,e;b,d)=\frac{(1-b)(1-d/qa)-(1-c)(1-e)\frac{bd}{qace}}{(1-b)(1-d/qa)},\\
&&\mathbb{A}_q:=\mathbb{A}_q(a,c,e;b,d)=\frac{(1-bd/qace)(1-qa)(1-c)(1-e)\,bd}{(1-b)(1-qb)(1-d)(1-d/qa)\,qace},\\
&&\mathfrak{A}_q:=\mathfrak{A}_q(a,c,e;b,d)=\frac{(1-bd/qace)}{(1-b)(1-d/qa)}.
 \enm
\end{thm}

Performing the substitutions $a\to q^a$, $b\to q^b$, $c\to q^c$,
$d\to q^d$, $e\to q^e$ for Theorem \ref{thm-a} and then letting
$q\to 1$, we recover the following relation.

\begin{corl}[\citu{wang}{Theorem 1}]\label{corl} For five complex numbers $\{a, b,
c, d, e\}$ subject to the condition $Re(b+d-a-c-e)>1$, there holds
the three-term contiguous relation of $_3F_2$-series:
 \bnm
&&_3F_2\ffnk{cccc}{1}{a,c,e}{b,d}
=\mathcal{A}\,{_3F_2}\ffnk{cccc}{1}{a+1,c,e}{b+1,d}
+\mathbb{A}\,{_3F_2}\ffnk{cccc}{1}{a+2,c+1,e+1}{b+2,d+1},\\
&&_3F_2\ffnk{cccc}{1}{a,c,e}{b,d}
=\mathcal{A}\,{_3F_2}\ffnk{cccc}{1}{a+1,c,e}{b+1,d}
+\mathfrak{A}\,{_3F_2}^{*}\ffnk{cccc}{1}{a+1,c,e}{b+1,d},
 \enm
where the coefficients $\mathcal{A}$, $\mathbb{A}$ and
$\mathfrak{A}$ are given by
 \bnm
&&\mathcal{A}:=\mathcal{A}(a,c,e;b,d)=\frac{(1+a-d)b+ce}{(1+a-d)b},\\
&&\mathbb{A}:=\mathbb{A}(a,c,e;b,d)=\frac{(1+a+c+e-b-d)(1+a)ce}{(1+a-d)(1+b)bd},\\
&&\mathfrak{A}:=\mathfrak{A}(a,c,e;b,d)=\frac{1+a+c+e-b-d}{(1+a-d)b}.
 \enm
\end{corl}
In Corollary \ref{corl}, the hypergeometric series and
\emph{shifted} hypergeometric series have been offered by
 \bnm
&&\xxqdn\qdn_{r}F_s\ffnk{cccc}{z}{a_1,&a_2,&\cdots,&a_r}{b_1,&b_2,&\cdots,&b_s}
 =\sum_{k=0}^\infty
\frac{(a_{1})_{k}(a_2)_{k}\cdots(a_r)_{k}}
 {(b_1)_{k}(b_2)_{k}\cdots(b_s)_{k}}\frac{z^k}{k!},\\
 &&\xxqdn\qdn\!_{r}F_s^{*}\ffnk{cccc}{z}{a_1,&a_2,&\cdots,&a_r}{b_1,&b_2,&\cdots,&b_s}
 =\sum_{k=0}^\infty
\frac{(a_{1})_{k}(a_2)_{k}\cdots(a_r)_{k}}
 {(b_1)_{k}(b_2)_{k}\cdots(b_s)_{k}}\frac{kz^k}{k!},
 \enm
where the shifted factorial is
 \[(x)_n=
\begin{cases}
\prod_{i=0}^{n-1}(x+i),&\quad\text{for}\quad n>0;\\
1,&\quad\text{for}\quad n=0;\\
\frac{1}{\prod_{j=n}^{-1}(x+j)},&\quad\text{for}\quad n<0.
\end{cases}\]

\subsection{Pattern B}$ $\\\\
For two sequences defined by
\[U_k=\ffnk{ccc}{q}{c,&e}{d/q,&qce/d}_k \quad\text{and}\quad V_k=\ffnk{ccc}{q}{a,&q^2ce/d}{q,&b}_k\bigg(\frac{bd}{qace}\bigg)^k,\]
we can easily verify the limiting relation
\[U_{0}V_{-1}=\lim_{n\to \infty}U_{n+1}V_{n}=0\]
and the finite differences
 \bnm
&&\xqdn\tilde{\triangle} U_k=\ffnk{ccc}{q}{c,&e}{d,&q^2ce/d}_kq^k\:\frac{(1-qc/d)(1-qe/d)}{(1-qce/d)(1-q/d)},\\
&&\xqdn\nabla V_k=\ffnk{ccc}{q}{a/q,&qce/d}{q,&b}_k
 \bigg(\frac{bd}{ace}\bigg)^k
 \bigg\{1+\frac{(1-qace/bd)}{(1-a/q)(1-qce/d)}\frac{1-q^k}{q^k}\bigg\}.
 \enm
By means of Lemma \ref{lemm}, we can reformulate the following
$_3\phi_2$-series:
 \bnm
&&\frac{(1-qc/d)(1-qe/d)}{(1-qce/d)(1-q/d)}{_3\phi_2}\ffnk{cccc}{q;\frac{bd}{ace}}{a,c,e}{b,d}
=\sum_{k\geq0}V_k\tilde{\triangle} U_k=\sum_{k\geq0}U_k\nabla V_k\\
&&=\:\sum_{k\geq0}\ffnk{ccc}{q}{a/q,&c,&e}{q,&b,&d/q}_k\bigg(\frac{bd}{ace}\bigg)^k
 \\&&+\:\frac{(1-qace/bd)}{(1-a/q)(1-qce/d)}
 \sum_{k\geq0}(1-q^k)\ffnk{ccc}{q}{a/q,&c,&e}{q,&b,&d/q}_k\bigg(\frac{bd}{qace}\bigg)^k.
 \enm
Shifting the summation index $k\to k+1$ for the last sum, we get the
following relation.

\begin{thm}[Pattern B]\label{thm-b} For five complex numbers $\{a, b,
c, d, e\}$ subject to the condition $|bd/qace|<1$, there holds the
three-term contiguous relation of $_3\phi_2$-series:
 \blm
&&_3\phi_2\ffnk{cccc}{q;\frac{bd}{ace}}{a,c,e}{b,d}
=\label{eq-b}\mathcal{B}_q{_3\phi_2}\ffnk{cccc}{q;\frac{bd}{ace}}{a/q,c,e}{b,d/q}
+\mathbb{B}_q{_3\phi_2}\ffnk{cccc}{q;\frac{bd}{qace}}{a,qc,qe}{qb,d},\\
&&_3\phi_2\ffnk{cccc}{q;\frac{bd}{ace}}{a,c,e}{b,d}
=\label{eq-bb}\mathcal{B}_q{_3\phi_2}\ffnk{cccc}{q;\frac{bd}{ace}}{a/q,c,e}{b,d/q}
+\mathfrak{B}_q{_3\phi_2}^{*}\ffnk{cccc}{q;\frac{bd}{qace}}{a/q,c,e}{b,d/q},
 \elm
where the coefficients $\mathcal{B}_q$, $\mathbb{B}_q$ and
$\mathfrak{B}_q$ are defined by
 \bnm
&&\mathcal{B}_q:=\mathcal{B}_q(a,c,e;b,d)=\frac{(1-qce/d)(1-q/d)}{(1-qc/d)(1-qe/d)},\\
&&\mathbb{B}_q:=\mathbb{B}_q(a,c,e;b,d)=\frac{(1-c)(1-e)(1-bd/qace)q}{(1-b)(1-qc/d)(1-qe/d)d},\\
&&\mathfrak{B}_q:=\mathfrak{B}_q(a,c,e;b,d)=\frac{(1-qace/bd)(1-q/d)}{(1-a/q)(1-qc/d)(1-qe/d)}.
 \enm
\end{thm}

Employing the substitutions $a\to q^a$, $b\to q^b$, $c\to q^c$,
$d\to q^d$, $e\to q^e$ for Theorem \ref{thm-b} and then letting
$q\to 1$, we recover the following relation.

\begin{corl}[\citu{wang}{Theorem 2}] For five complex numbers $\{a, b,
c, d, e\}$ subject to the condition $Re(b+d-a-c-e)>1$, there holds
the three-term contiguous relation of $_3F_2$-series:
 \bnm
&&_3F_2\ffnk{cccc}{1}{a,c,e}{b,d}
=\mathcal{B}\,{_3F_2}\ffnk{cccc}{1}{a-1,c,e}{b,d-1}
+\mathbb{B}\,{_3F_2}\ffnk{cccc}{1}{a,c+1,e+1}{b+1,d},\\
&&_3F_2\ffnk{cccc}{1}{a,c,e}{b,d}
=\mathcal{B}\,{_3F_2}\ffnk{cccc}{1}{a-1,c,e}{b,d-1}
+\mathfrak{B}\,{_3F_2}^{*}\ffnk{cccc}{1}{a-1,c,e}{b,d-1},
 \enm
where the coefficients $\mathcal{B}$, $\mathbb{B}$ and
$\mathfrak{B}$ are given by
 \bnm
&&\mathcal{B}:=\mathcal{B}(a,c,e;b,d)=\frac{(1+c+e-d)(1-d)}{(1+c-d)(1+e-d)},\\
&&\mathbb{B}:=\mathbb{B}(a,c,e;b,d)=\frac{(1+a+c+e-b-d)ce}{(1+c-d)(d-e-1)b},\\
&&\mathfrak{B}:=\mathfrak{B}(a,c,e;b,d)=\frac{(1+a+c+e-b-d)(1-d)}{(1-a)(1+c-d)(d-e-1)}.
 \enm
\end{corl}
\subsection{Pattern C}$ $\\\\
Define two sequences by
\[U_k=\ffnk{ccc}{q}{qc,&qe}{q,&qce}_k \quad\text{and}\quad
 V_k=\ffnk{ccc}{q}{a,&qce}{b,&d}_k\bigg(\frac{bd}{qace}\bigg)^k.\]
Then it is not hard to check the limiting relation
\[U_{-1}V_{0}=\lim_{n\to \infty}U_{n}V_{n+1}=0\]
and the finite differences
 \bnm
&&\xqdn\nabla U_k=\ffnk{ccc}{q}{c,&e}{q,&qce}_kq^k,\\
&&\xqdn\tilde{\triangle} V_k=\ffnk{ccc}{q}{a,&qce}{qb,&qd}_k
 \bigg(\frac{bd}{ace}\bigg)^k
 \bigg\{\frac{\sss(1-b)(1-d)-(1-a)(1-qce)\frac{bd}{qace}}{(1-b)(1-d)}
 +\frac{(1-bd/qace)}{(1-b)(1-d)}\frac{1-q^k}{q^k}\bigg\}.
 \enm
According to Lemma \ref{lemm}, we can recombine the following
$_3\phi_2$-series:
 \bnm
&&_3\phi_2\ffnk{cccc}{q;\frac{bd}{ace}}{a,c,e}{b,d}
=\sum_{k\geq0}V_k\nabla U_k=\sum_{k\geq0}U_k\tilde{\triangle} V_k\\
&&\:\:=\:\:\frac{\sss(1-b)(1-d)-(1-a)(1-qce)\frac{bd}{qace}}{(1-b)(1-d)}
\sum_{k\geq0}\ffnk{ccc}{q}{a,&qc,&qe}{q,&qb,&qd}_k\bigg(\frac{bd}{ace}\bigg)^k\\
 &&\:\:+\:\:\:\frac{(1-bd/qace)}{(1-b)(1-d)}
 \sum_{k\geq0}(1-q^k)\ffnk{ccc}{q}{a,&qc,&qe}{q,&qb,&qd}_k
 \bigg(\frac{bd}{qace}\bigg)^k.
 \enm
Shifting the summation index $k\to k+1$ for the last sum, we derive
the following relation.

\begin{thm}[Pattern C]\label{thm-c} For five complex numbers $\{a, b,
c, d, e\}$ subject to the condition $|bd/qace|<1$, there holds the
three-term contiguous relation of $_3\phi_2$-series:
 \blm
&&_3\phi_2\ffnk{cccc}{q;\frac{bd}{ace}}{a,c,e}{b,d}
\label{eq-c}=\mathcal{C}_q{_3\phi_2}\ffnk{cccc}{q;\frac{bd}{ace}}{a,qc,qe}{qb,qd}
+\mathbb{C}_q{_3\phi_2}\ffnk{cccc}{q;\frac{bd}{qace}}{qa,q^2c,q^2e}{q^2b,q^2d},\\
&&_3\phi_2\ffnk{cccc}{q;\frac{bd}{ace}}{a,c,e}{b,d}
=\label{eq-cc}\mathcal{C}_q{_3\phi_2}\ffnk{cccc}{q;\frac{bd}{ace}}{a,qc,qe}{qb,qd}
+\mathfrak{C}_q{_3\phi_2}^{*}\ffnk{cccc}{q;\frac{bd}{qace}}{a,qc,qe}{qb,qd},
 \elm
where the coefficients $\mathcal{C}_q$, $\mathbb{C}_q$ and
$\mathfrak{C}_q$ are defined by
 \bnm
&&\mathcal{C}_q:=\mathcal{C}_q(a,c,e;b,d)=\frac{(1-b)(1-d)-(1-a)(1-qce)\frac{bd}{qace}}{(1-b)(1-d)},\\
&&\mathbb{C}_q:=\mathbb{C}_q(a,c,e;b,d)=\frac{(1-bd/qace)(1-a)(1-qc)(1-qe)\,bd}{(1-b)(1-d)(1-qb)(1-qd)\,qace},\\
&&\mathfrak{C}_q:=\mathfrak{C}_q(a,c,e;b,d)=\frac{(1-bd/qace)}{(1-b)(1-d)}.
 \enm
\end{thm}

Performing the substitutions $a\to q^a$, $b\to q^b$, $c\to q^c$,
$d\to q^d$, $e\to q^e$ for Theorem \ref{thm-c} and then letting
$q\to 1$, we recover the following relation.

\begin{corl}[\citu{wang}{Theorem 3}] For five complex numbers $\{a, b,
c, d, e\}$ subject to the condition $Re(b+d-a-c-e)>1$, there holds
the three-term contiguous relation of $_3F_2$-series:
 \bnm
&&_3F_2\ffnk{cccc}{1}{a,c,e}{b,d}
=\mathcal{C}\,{_3F_2}\ffnk{cccc}{1}{a,c+1,e+1}{b+1,d+1}
+\mathbb{C}\,{_3F_2}\ffnk{cccc}{1}{a+1,c+2,e+2}{b+2,d+2},\\
&&_3F_2\ffnk{cccc}{1}{a,c,e}{b,d}
=\mathcal{C}\,{_3F_2}\ffnk{cccc}{1}{a,c+1,e+1}{b+1,d+1}
+\mathfrak{C}\,{_3F_2}^{*}\ffnk{cccc}{1}{a,c+1,e+1}{b+1,d+1},
 \enm
where the coefficients $\mathcal{C}$, $\mathbb{C}$ and
$\mathfrak{C}$ are given by
 \bnm
&&\mathcal{C}:=\mathcal{C}(a,c,e;b,d)=\frac{bd-a(1+c+e)}{bd},\\
&&\mathbb{C}:=\mathbb{C}(a,c,e;b,d)=\frac{(b+d-a-c-e-1)(1+c)(1+e)a}{(1+b)(1+d)bd},\\
&&\mathfrak{C}:=\mathfrak{C}(a,c,e;b,d)=\frac{b+d-a-c-e-1}{bd}.
 \enm
\end{corl}

\subsection{Pattern D}$ $\\\\
For two sequences defined by
\[U_k=\ffnk{ccc}{q}{a,&bd/q^2a}{b/q,&d/q}_k \quad\text{and}\quad V_k=\ffnk{ccc}{q}{c,&e}{q,&bd/q^2a}_k\bigg(\frac{bd}{qace}\bigg)^k,\]
we can verify without difficulty the limiting relation
\[U_{0}V_{-1}=\lim_{n\to \infty}U_{n+1}V_{n}=0\]
and the finite differences
 \bnm
&&\xqdn\tilde{\triangle} U_k=\ffnk{ccc}{q}{a,&bd/q^2a}{b,&d}_kq^k\:\frac{(1-qa/b)(1-qa/d)}{(1-q/b)(1-q/d)a},\\
&&\xqdn\nabla V_k=\ffnk{ccc}{q}{c/q,&e/q}{q,&bd/q^2a}_k
 \bigg(\frac{bd}{ace}\bigg)^k
 \bigg\{1+\frac{(1-qace/bd)}{(1-c/q)(1-e/q)}\frac{1-q^k}{q^k}\bigg\}.
 \enm
In terms of Lemma \ref{lemm}, we can recompose the following
$_3\phi_2$-series:
 \bnm
&&\frac{(1-qa/b)(1-qa/d)}{(1-q/b)(1-q/d)a}{_3\phi_2}\ffnk{cccc}{q;\frac{bd}{ace}}{a,c,e}{b,d}
=\sum_{k\geq0}V_k\tilde{\triangle} U_k=\sum_{k\geq0}U_k\nabla V_k\\
&&=\:\sum_{k\geq0}\ffnk{ccc}{q}{a,&c/q,&e/q}{q,&b/q,&d/q}_k\bigg(\frac{bd}{ace}\bigg)^k
 \\&&+\:\frac{(1-qace/bd)}{(1-c/q)(1-e/q)}
 \sum_{k\geq0}(1-q^k)\ffnk{ccc}{q}{a,&c/q,&e/q}{q,&b/q,&d/q}_k\bigg(\frac{bd}{qace}\bigg)^k.
 \enm
Shifting the summation index $k\to k+1$ for the last sum, we deduce
the following relation.

\begin{thm}[Pattern D]\label{thm-d} For five complex numbers $\{a, b,
c, d, e\}$ subject to the condition $|bd/qace|<1$, there holds the
three-term contiguous relation of $_3\phi_2$-series:
 \blm
&&\qqdn_3\phi_2\ffnk{cccc}{q;\frac{bd}{ace}}{a,c,e}{b,d}
=\label{eq-d}\mathcal{D}_q{_3\phi_2}\ffnk{cccc}{q;\frac{bd}{ace}}{a,c/q,e/q}{b/q,d/q}
+\mathbb{D}_q{_3\phi_2}\ffnk{cccc}{q;\frac{bd}{qace}}{qa,c,e}{b,d},\\
&&\qqdn_3\phi_2\ffnk{cccc}{q;\frac{bd}{ace}}{a,c,e}{b,d}
=\label{eq-dd}\mathcal{D}_q{_3\phi_2}\ffnk{cccc}{q;\frac{bd}{ace}}{a,c/q,e/q}{b/q,d/q}
+\mathfrak{D}_q{_3\phi_2}^{*}\ffnk{cccc}{q;\frac{bd}{qace}}{a,c/q,e/q}{b/q,d/q},
 \elm
where the coefficients $\mathcal{D}_q$, $\mathbb{D}_q$ and
$\mathfrak{D}_q$ are defined by
 \bnm
&&\mathcal{D}_q:=\mathcal{D}_q(a,c,e;b,d)=\frac{(1-q/b)(1-q/d)a}{(1-qa/b)(1-qa/d)},\\
&&\mathbb{D}_q:=\mathbb{D}_q(a,c,e;b,d)=\frac{(1-a)(1-qace/bd)q}{(1-qa/b)(1-qa/d)ce},\\
&&\mathfrak{D}_q:=\mathfrak{D}_q(a,c,e;b,d)=\frac{(1-qace/bd)(1-q/b)(1-q/d)a}{(1-qa/b)(1-qa/d)(1-c/q)(1-e/q)}.
 \enm
\end{thm}

Employing the substitutions $a\to q^a$, $b\to q^b$, $c\to q^c$,
$d\to q^d$, $e\to q^e$ for Theorem \ref{thm-d} and then letting
$q\to 1$, we recover the following relation.

\begin{corl}[\citu{wang}{Theorem 4}] For five complex numbers $\{a, b,
c, d, e\}$ subject to the condition $Re(b+d-a-c-e)>1$, there holds
the three-term contiguous relation of $_3F_2$-series:
 \bnm
&&_3F_2\ffnk{cccc}{1}{a,c,e}{b,d}
=\mathcal{D}\,{_3F_2}\ffnk{cccc}{1}{a,c-1,e-1}{b-1,d-1}
+\mathbb{D}\,{_3F_2}\ffnk{cccc}{1}{a+1,c,e}{b,d},\\
&&_3F_2\ffnk{cccc}{1}{a,c,e}{b,d}
=\mathcal{D}\,{_3F_2}\ffnk{cccc}{1}{a,c-1,e-1}{b-1,d-1}
+\mathfrak{D}\,{_3F_2}^{*}\ffnk{cccc}{1}{a,c-1,e-1}{b-1,d-1},
 \enm
where the coefficients $\mathcal{D}$, $\mathbb{D}$ and
$\mathfrak{D}$ are given by
 \bnm
&&\mathcal{D}:=\mathcal{D}(a,c,e;b,d)=\frac{(1-b)(1-d)}{(1+a-b)(1+a-d)},\\
&&\mathbb{D}:=\mathbb{D}(a,c,e;b,d)=\frac{a(1+a+c+e-b-d)}{(1+a-b)(1+a-d)},\\
&&\mathfrak{D}:=\mathfrak{D}(a,c,e;b,d)=\frac{(1+a+c+e-b-d)(1-b)(1-d)}{(1+a-b)(1+a-d)(1-c)(1-e)}.
 \enm
\end{corl}

\textbf{Remark:} The convergent conditions for the contiguous
relations that will emergence in the next two sections are easy to
be confirmed. For simplifying the expressions, we shall not lay out
them one by one.
\section{Nine two-term contiguous relations of $_3\phi_2$-series}
In terms of one or two patterns of \textbf{A}, \textbf{B},
\textbf{C} and \textbf{D}, we offer nine two-term contiguous
relations of $_3\phi_2$-series in this section. Each subsection will
be labeled by the corresponding patterns.

\subsection{\textbf{A}\&\textbf{A}} $ $\\\\
Let Eq$^{\star}$\eqref{eq-a} stand for Eq\eqref{eq-a} under the
parameter replacements
\[a\to c/q,\quad c\to qa,\quad b\to d/q,\quad d\to qb.\]
Then consider the linear combination of two equations
\[\text{Eq}\eqref{eq-a}-\text{Eq}^{\star}\eqref{eq-a}\frac{(1-qb/c)(1-d/q)}{(1-b)(1-d/qa)}
\quad\text{with}\quad d=\frac{qabe(q-c)}{qab+ce-bc-ace}.\]
 With this specific value $d$, we can check that the right member of
 the last equation vanishes. After some simplification, we attain
 the following relation.

\begin{thm}[Two-term contiguous relation of $_3\phi_2$-series]\label{thm-f}
 \bnm
&&_3\phi_2\ffnk{cccc}{q;\frac{qb^2(q-c)}{c(qab+ce-bc-ace)}}{a,c,e}{b,\frac{qabe(q-c)}{qab+ce-bc-ace}}\\
&=&_3\phi_2\ffnk{cccc}{q;\frac{qb^2(q-c)}{c(qab+ce-bc-ace)}}{qa,c/q,e}{qb,\frac{abe(q-c)}{qab+ce-bc-ace}}\\
&\times&\frac{(1-qb/c)(qab+abce+ce-bc-qabe-ace)}{(1-b)(qab+bce+ce-bc-qbe-ace)}.
 \enm
\end{thm}

Performing the substitutions $a\to q^a$, $b\to q^b$, $c\to q^c$,
$e\to q^e$ for Theorem \ref{thm-f} and then letting $q\to 1$, we
recover the following relation.

\begin{corl}[\citu{wang}{Theorem 5}]
 \bnm
_3F_2\ffnk{cccc}{1}{a,c,e}{b,1+e-\frac{a(b-e)}{1-c}}
&=&{_3F_2}\ffnk{cccc}{1}{a+1,c-1,e}{b+1,e-\frac{a(b-e)}{1-c}}\\
&\times&\frac{(1+b-c)(ab+ce-ae-e)}{b(a+ab+ce-ac-ae-e)}.
 \enm
\end{corl}

Specifying the parameters, in Theorem \ref{thm-f},  by
\[a\to \alpha,\quad c\to \beta,\quad e\to\gamma,\quad b\to q\alpha,\]
we achieve the following relation with one free parameter less.

\begin{prop}[Two-term contiguous relation of $_3\phi_2$-series]\label{prop-a}
 \bnm
&&_3\phi_2\ffnk{cccc}{q;\frac{q^3\alpha^2(q-\beta)}{\beta(q^2\alpha^2+\beta\gamma-q\alpha\beta-\alpha\beta\gamma)}}
{\alpha,\beta,\gamma}{q\alpha,\frac{q^2\alpha^2\gamma(q-\beta)}{q^2\alpha^2+\beta\gamma-q\alpha\beta-\alpha\beta\gamma}}\\
&=&_3\phi_2\ffnk{cccc}{q;\frac{q^3\alpha^2(q-\beta)}{\beta(q^2\alpha^2+\beta\gamma-q\alpha\beta-\alpha\beta\gamma)}}
{q\alpha,\beta/q,\gamma}{q^2\alpha,\frac{q\alpha^2\gamma(q-\beta)}{q^2\alpha^2+\beta\gamma-q\alpha\beta-\alpha\beta\gamma}}\\
&\times&\frac{(1-q^2\alpha/\beta)(q^2\alpha^2+q\alpha^2\beta\gamma+\beta\gamma-q^2\alpha^2\gamma-q\alpha\beta-\alpha\beta\gamma)}
{(1-q\alpha)(q^2\alpha^2+q\alpha\beta\gamma+\beta\gamma-q^2\alpha\gamma-q\alpha\beta-\alpha\beta\gamma)}.
 \enm
\end{prop}

Employing the substitutions $\alpha\to q^\alpha$, $\beta\to
q^\beta$, $\gamma\to q^\gamma$ for Proposition \ref{prop-a} and then
letting $q\to 1$, we recover the following relation.

\begin{corl}[\citu{krattenthaler}{Theorem 2}, see also \citu{wang}{Proposition 6}]
 \bnm
_3F_2\ffnk{cccc}{1}{\alpha,\beta,\gamma}{\alpha+1,\gamma+\frac{\alpha(\alpha-\gamma+1)}{\beta-1}+1}
&=&{_3F_2}\ffnk{cccc}{1}{\alpha+1,\beta-1,\gamma}{\alpha+2,\gamma+\frac{\alpha(\alpha-\gamma+1)}{\beta-1}}\\
&\times&\frac{(2+\alpha-\beta)(\alpha+\alpha^2+\beta\gamma-\alpha\gamma-\gamma)}
{(1+\alpha)(2\alpha+\alpha^2+\beta\gamma-\alpha\beta-\alpha\gamma-\gamma)}.
 \enm
\end{corl}
\subsection{\textbf{A}\&\textbf{B}}$ $\\\\
Let Eq$^{\star}$\eqref{eq-b} stand for Eq\eqref{eq-b} under the
parameter replacements
\[a\to q^2a,\quad b\to d,\quad d\to q^2b.\]
Then consider the linear combination of two equations
\[\text{Eq}\eqref{eq-a}-\text{Eq}^{\star}\eqref{eq-b}\frac{(1-qa)(1-qb/c)(1-qb/e)}{(1-b)(1-qb)(qa/d-1)q}
\quad\text{with}\quad
d=\frac{q^2ace(1-b)}{q^2ab+qce+ce-qace-qbc-qbe}.\]
 With this specific value $d$, we can verify that the right member of
 the last equation vanishes. After some simplification, we establish
 the following relation.

\begin{thm}[Two-term contiguous relation of $_3\phi_2$-series]\label{thm-g}
 \bnm
&&_3\phi_2\ffnk{cccc}{q;\frac{q^2b(1-b)}{q^2ab+qce+ce-qace-qbc-qbe}}{a,c,e}
{b,\quad\frac{q^2ace(1-b)}{q^2ab+qce+ce-qace-qbc-qbe}}\\
&=&_3\phi_2\ffnk{cccc}{q;\frac{q^2b(1-b)}{q^2ab+qce+ce-qace-qbc-qbe}}{q^2a,c,e}
{q^2b,\frac{q^2ace(1-b)}{q^2ab+qce+ce-qace-qbc-qbe}}\\
&\times&\frac{(1-qa)(c-qb)(e-qb)}{(1-qb)(q^2ab+qbce+ce-qace-qbc-qbe)}.
 \enm
\end{thm}

Performing the substitutions $a\to q^a$, $b\to q^b$, $c\to q^c$,
$e\to q^e$ for Theorem \ref{thm-g} and then letting $q\to 1$, we
recover the following relation.

\begin{corl}[\citu{wang}{Theorem 7}]
 \bnm
_3F_2\ffnk{cccc}{1}{a,c,e}{b,2+2a+\frac{(a+1)(1-c-e)+ce}{b}}
&=&{_3F_2}\ffnk{cccc}{1}{a+2,c,e}{b+2,2+2a+\frac{(a+1)(1-c-e)+ce}{b}}\\
&\times&\frac{(1+a)(1+b-c)(1+b-e)}{(1+b)(1+a+b+ab+ce-c-e-ac-ae)}.
 \enm
\end{corl}

Specifying the parameters, in Theorem \ref{thm-g},  by
\[a\to \alpha,\quad c\to \beta,\quad e\to\gamma,\quad b\to q\alpha,\]
we found the following relation with one free parameter less.

\begin{prop}[Two-term contiguous relation of $_3\phi_2$-series]\label{prop-b}
 \bnm
&&\qdn_3\phi_2\ffnk{cccc}{q;\frac{q^3\alpha(1-q\alpha)}{q^3\alpha^2+q\beta\gamma+\beta\gamma-q\alpha\beta\gamma-q^2\alpha\beta-q^2\alpha\gamma}}
{\alpha,\beta,\gamma}{q\alpha,\:\:\frac{q^2\alpha\beta\gamma(1-q\alpha)}{q^3\alpha^2+q\beta\gamma+\beta\gamma-q\alpha\beta\gamma-q^2\alpha\beta-q^2\alpha\gamma}}\\
&&\xqdn=\:{_3\phi_2}\ffnk{cccc}{q;\frac{q^3\alpha(1-q\alpha)}{q^3\alpha^2+q\beta\gamma+\beta\gamma-q\alpha\beta\gamma-q^2\alpha\beta-q^2\alpha\gamma}}
{q^2\alpha,\beta,\gamma}{q^3\alpha,\frac{q^2\alpha\beta\gamma(1-q\alpha)}{q^3\alpha^2+q\beta\gamma+\beta\gamma-q\alpha\beta\gamma-q^2\alpha\beta-q^2\alpha\gamma}}\\
&&\xqdn\times\:\frac{(1-q\alpha)(\beta-q^2\alpha)(\gamma-q^2\alpha)}
{(1-q^2\alpha)(q^3\alpha^2+q^2\alpha\beta\gamma+\beta\gamma-q\alpha\beta\gamma-q^2\alpha\beta-q^2\alpha\gamma)}.
 \enm
\end{prop}

Employing the substitutions $\alpha\to q^\alpha$, $\beta\to
q^\beta$, $\gamma\to q^\gamma$ for Proposition \ref{prop-b} and then
letting $q\to 1$, we recover the following relation.

\begin{corl}[\citu{krattenthaler}{Theorem 10}, see also \citu{wang}{Proposition 8}]
 \bnm
_3F_2\ffnk{cccc}{1}{\alpha,\beta,\gamma}{\alpha+1,3+2\alpha-\beta-\gamma+\frac{\beta\gamma}{\alpha+1}}
&=&_3F_2\ffnk{cccc}{1}{\alpha+2,\beta,\gamma}{\alpha+3,3+2\alpha-\beta-\gamma+\frac{\beta\gamma}{\alpha+1}}\\
&\times&\frac{(1+\alpha)(\alpha-\beta+2)(\alpha-\gamma+2)}
{(2+\alpha)(2+3\alpha+\alpha^2+\beta\gamma-\alpha\beta-\alpha\gamma-\beta-\gamma)}.
 \enm
\end{corl}

\subsection{\textbf{A}\&\textbf{B}}$ $\\\\
Let Eq$^{\star}$\eqref{eq-b} stand for Eq\eqref{eq-b} under the
parameter replacements
\[a\to q^2a,\quad b\to qb,\quad d\to qd.\]
Then consider the linear combination of two equations
\[\text{Eq}\eqref{eq-a}-\text{Eq}^{\star}\eqref{eq-b}\frac{(1-qa)(1-d/c)(1-d/e)b}{(1-b)(1-d)(qa-d)}
\quad\text{with}\quad d=\frac{ce(qa-b)}{qab+ce-bc-be}.\]
 With this specific value $d$, we can check that the right member of
 the last equation vanishes. After some simplification, we obtain
 the following relation.

\begin{thm}[Two-term contiguous relation of $_3\phi_2$-series]\label{thm-h}
 \bnm
&&_3\phi_2\ffnk{cccc}{q;\frac{b(qa-b)}{a(qab+ce-bc-be)}}{a,c,e}
{b,\:\:\frac{ce(qa-b)}{qab+ce-bc-be}}\\
&=&_3\phi_2\ffnk{cccc}{q;\frac{b(qa-b)}{a(qab+ce-bc-be)}}{q^2a,c,e}
{qb,\frac{qce(qa-b)}{qab+ce-bc-be}}\\
&\times&\frac{(1-qa)(b-c)(b-e)}{(1-b)(qab+bce+ce-bc-be-qace)}.
 \enm
\end{thm}

Performing the substitutions $a\to q^a$, $b\to q^b$, $c\to q^c$,
$e\to q^e$ for Theorem \ref{thm-h} and then letting $q\to 1$, we
recover the following relation.

\begin{corl}[\citu{wang}{Theorem 9}]
 \bnm
_3F_2\ffnk{cccc}{1}{a,c,e}{b,1+a-\frac{(1+a-c)(1+a-e)}{1+a-b}}
&=&{_3F_2}\ffnk{cccc}{1}{a+2,c,e}{b+1,2+a-\frac{(1+a-c)(1+a-e)}{1+a-b}}\\
&\times&\frac{(1+a)(b-c)(b-e)}{b(ab+ce+b-c-e-ac-ae)}.
 \enm
\end{corl}

\subsection{\textbf{A}\&\textbf{B}}$ $\\\\
Let Eq$^{\star}$\eqref{eq-b} stand for Eq\eqref{eq-b} under the
parameter replacements
\[a\to qc,\quad c\to qa,\quad b\to d,\quad d\to q^2b.\]
Then consider the linear combination of two equations
\[\text{Eq}\eqref{eq-a}-\text{Eq}^{\star}\eqref{eq-b}\frac{(b-a)(1-c)(1-qb/e)}{(1-b)(1-qb)(1-qa/d)c}
\quad\text{with}\quad d=\frac{qac(1-b)}{a+c-ac-b}.\]
 With this specific value $d$, we can verify that the right member of
 the last equation vanishes. After some simplification, we get
 the following relation.

\begin{thm}[Two-term contiguous relation of $_3\phi_2$-series]\label{thm-i}
 \bnm
\:_3\phi_2\ffnk{cccc}{q;\frac{qb(1-b)}{(a+c-ac-b)e}}{a,c,e}
{b,\:\:\frac{qac(1-b)}{a+c-ac-b}}=
{_3\phi_2}\ffnk{cccc}{q;\frac{qb(1-b)}{(a+c-ac-b)e}}{qa,qc,e}
{q^2b,\frac{qac(1-b)}{a+c-ac-b}}\frac{qb-e}{qbe-e}.
 \enm
\end{thm}

Employing the substitutions $a\to q^a$, $b\to q^b$, $c\to q^c$,
$e\to q^e$ for Theorem \ref{thm-i} and then letting $q\to 1$, we
recover the following relation.

\begin{corl}[\citu{wang}{Theorem 10}]
 \bnm
_3F_2\ffnk{cccc}{1}{a,c,e}{b,1+a+c-\frac{ac}{b}}
={_3F_2}\ffnk{cccc}{1}{a+1,c+1,e}{b+2,1+a+c-\frac{ac}{b}}\frac{1+b-e}{1+b}.
 \enm
\end{corl}

Specifying the parameters, in Theorem \ref{thm-i},  by
\[a\to \beta,\quad c\to \gamma,\quad e\to \beta\gamma/\alpha,\quad b\to \frac{\beta+\gamma-\alpha-\beta\gamma}{1-\alpha},\]
we derive the following relation with one free parameter less.

\begin{prop}[Two-term contiguous relation of $_3\phi_2$-series]\label{prop-c}
 \bnm
&&_3\phi_2\ffnk{cccc}{q;\frac{q(\beta+\gamma-\alpha-\beta\gamma)}{(1-\alpha)\beta\gamma}}
{\beta\gamma/\alpha,\beta,\gamma}{\frac{q\beta\gamma}{\alpha},\quad\:\:\frac{\beta+\gamma-\alpha-\beta\gamma}{1-\alpha}}\\
&=&_3\phi_2\ffnk{cccc}{q;\frac{q(\beta+\gamma-\alpha-\beta\gamma)}{(1-\alpha)\beta\gamma}}
{\beta\gamma/\alpha,q\beta,q\gamma}{\frac{q\beta\gamma}{\alpha},\frac{q^2(\beta+\gamma-\alpha-\beta\gamma)}{1-\alpha}}\\
&\times&\frac{\alpha(q\alpha+q\beta\gamma-\beta\gamma-q\beta-q\gamma)+\beta\gamma}
{\alpha(1+q\alpha+q\beta\gamma-\alpha-q\beta-q\gamma)}.
 \enm
\end{prop}

Performing the substitutions $\alpha\to q^\alpha$, $\beta\to
q^\beta$, $\gamma\to q^\gamma$  for Proposition \ref{prop-c} and
then letting $q\to1$, we recover the following relation.

\begin{corl}[\citu{krattenthaler}{Theorem 9}, see also \citu{wang}{Proposition 11}]
 \bnm
\qquad
_3F_2\ffnk{cccc}{1}{\beta+\gamma-\alpha,\beta,\gamma}{\beta+\gamma-\alpha+1,\frac{\beta\gamma}{\alpha}}
={_3F_2}\ffnk{cccc}{1}{\beta+\gamma-\alpha,\beta+1,\gamma+1}{\beta+\gamma-\alpha+1,\frac{\beta\gamma}{\alpha}+2}
\frac{\alpha+\alpha^2+\beta\gamma-\alpha\beta-\alpha\gamma}
{\alpha+\beta\gamma}.
 \enm
\end{corl}

\subsection{\textbf{A}\&\textbf{C}} $ $\\\\
Let Eq$^{\star}$\eqref{eq-c} stand for Eq\eqref{eq-c} under the
parameter replacements
\[a\to qa,\quad c\to c/q,\quad e\to e/q,\quad d\to d/q.\]
Then consider the linear combination of two equations
\[\text{Eq}\eqref{eq-a}-\text{Eq}^{\star}\eqref{eq-c}\frac{1-d/q}{1-d/qa}
\quad\text{with}\quad b=\frac{qce(a-1)}{q^2a+ce-qc-qe}.\]
 With this specific value $b$, we can check that the right member of
 the last equation vanishes. After some simplification, we deduce
 the following relation.

\begin{thm}[Two-term contiguous relation of $_3\phi_2$-series]\label{thm-j}
 \bnm
&&_3\phi_2\ffnk{cccc}{q;\frac{qd(a-1)}{a(q^2a+ce-qc-qe)}}{a,c,e}
{d,\:\:\frac{qce(a-1)}{q^2a+ce-qc-qe}}\\
&=&_3\phi_2\ffnk{cccc}{q;\frac{qd(a-1)}{a(q^2a+ce-qc-qe)}}{qa,c/q,e/q}
{d/q,\frac{qce(a-1)}{q^2a+ce-qc-qe}}\frac{1-d/q}{1-d/qa}.
 \enm
\end{thm}

Employing the substitutions $a\to q^a$, $c\to q^c$, $d\to q^d$,
$e\to q^e$ for Theorem \ref{thm-j} and then letting $q\to1$, we
recover the following relation.

\begin{corl}[\citu{wang}{Theorem 12}]
 \bnm
\quad_3F_2\ffnk{cccc}{1}{a,c,e}{d,c+e-1-\frac{(1-c)(1-e)}{a}}
={_3F_2}\ffnk{cccc}{1}{a+1,c-1,e-1}{d-1,c+e-1-\frac{(1-c)(1-e)}{a}}\frac{1-d}{1+a-d}.
 \enm
\end{corl}

\subsection{\textbf{A}\&\textbf{D}}$ $\\\\
Let Eq$^{\star}$\eqref{eq-d} stand for Eq\eqref{eq-d} under the
parameter replacements
\[a\to c,\quad c\to q^2a,\quad e\to qe,\quad b\to q^2b,\quad d\to qd.\]
Then consider the linear combination of two equations
\[\text{Eq}\eqref{eq-a}-\text{Eq}^{\star}\eqref{eq-d}\frac{(1-qa)(1-qb/c)(1-d/c)(1-e)b}{(1-b)(1-d)(1-qb)(1-qa/d)e}
\quad\text{with}\quad d=\frac{qace(b-1)}{qabe+bc-qab-ce}.\]
 With this specific value $d$, we can verify that the right member of
 the last equation vanishes. After some simplification, we attain
 the following relation.

\begin{thm}[Two-term contiguous relation of $_3\phi_2$-series]\label{thm-k}
 \bnm
&&_3\phi_2\ffnk{cccc}{q;\frac{qb(b-1)}{qabe+bc-qab-ce}}{a,c,e}
{b,\:\:\frac{qace(b-1)}{qabe+bc-qab-ce}}\\
&=&_3\phi_2\ffnk{cccc}{q;\frac{qb(b-1)}{qabe+bc-qab-ce}}{q^2a,c,qe}
{q^2b,\frac{q^2ace(b-1)}{qabe+bc-qab-ce}}\\
&\times&\frac{(1-qa)(b-e)(c-qb)}{(1-qb)(qabe+qace+bc-qabce-qab-ce)}.
 \enm
\end{thm}

Performing the substitutions $a\to q^a$, $b\to q^b$, $c\to q^c$,
$e\to q^e$ for Theorem \ref{thm-k} and then letting $q\to1$, we
recover the following relation.

\begin{corl}[\citu{wang}{Theorem 13}]
 \bnm
_3F_2\ffnk{cccc}{1}{a,c,e}{b,1+a+\frac{e(c-a-1)}{b}}
&=&{_3F_2}\ffnk{cccc}{1}{a+2,c,e+1}{b+2,2+a+\frac{e(c-a-1)}{b}}\\
&\times&\frac{(1+a)(b-e)(1+b-c)}{(1+b)(ab+ce+b-ae-e)}.
 \enm
\end{corl}

\subsection{\textbf{B}\&\textbf{C}}$ $\\\\
Let Eq$^{\star}$\eqref{eq-c} stand for Eq\eqref{eq-c} under the
parameter replacements
\[a\to c,\quad c\to a/q^2,\quad e\to e/q,\quad b\to b/q,\quad d\to d/q^2.\]
Then consider the linear combination of two equations
\[\text{Eq}\eqref{eq-b}-\text{Eq}^{\star}\eqref{eq-c}\frac{(1-q/b)(1-q/d)(1-q^2/d)ce}{(1-q/a)(1-qc/d)(1-qe/d)q}
\quad\text{with}\quad d=\frac{q^3ce(a-b)}{qab+qace-q^2bc-abe}.\]
 With this specific value $d$, we can check that the right member of
 the last equation vanishes. After some simplification, we achieve
 the following relation.

\begin{thm}[Two-term contiguous relation of $_3\phi_2$-series]\label{thm-m}
 \bnm
&&_3\phi_2\ffnk{cccc}{q;\frac{q^3b(a-b)}{a(qab+qace-q^2bc-abe)}}{a,c,e}
{b,\:\:\frac{q^3ce(a-b)}{qab+qace-q^2bc-abe}}\\
&=&_3\phi_2\ffnk{cccc}{q;\frac{q^3b(a-b)}{a(qab+qace-q^2bc-abe)}}{a/q^2,c,e/q}
{b/q,\frac{qce(a-b)}{qab+qace-q^2bc-abe}}\\
&\times&\frac{(a-qc)(q-b)(qab+qace+q^2bce-q^2ace-q^2bc-abe)}
{(q-a)(b-qc)(qab+qace+q^2be-q^2ae-q^2bc-abe)}.
 \enm
\end{thm}

Employing the substitutions $a\to q^a$, $b\to q^b$, $c\to q^c$,
$e\to q^e$ for Theorem \ref{thm-m} and then letting $q\to1$, we
recover the following relation.

\begin{corl}[\citu{wang}{Theorem 14}]
 \bnm
\qquad_3F_2\ffnk{cccc}{1}{a,c,e}{b,2+\frac{(1-e)(1+c-a)}{a-b}}
={_3F_2}\ffnk{cccc}{1}{a-2,c,e-1}{b-1,\frac{(1-e)(1+c-a)}{a-b}}\\
\qquad\times\frac{(1-b)(1+c-a)(1+c+ae-ce-b-e)}{(1-a)(1+c-b)(1+c+ae+bc-ac-ce-b-e)}.
 \enm
\end{corl}

\subsection{\textbf{C}\&\textbf{D}}$ $\\\\
Let Eq$^{\star}$\eqref{eq-d} stand for Eq\eqref{eq-d} under the
parameter replacements
\[c\to q^2c,\quad e\to q^2e,\quad b\to q^2b,\quad d\to q^2d.\]
Then consider the linear combination of two equations
\[\text{Eq}\eqref{eq-c}-\text{Eq}^{\star}\eqref{eq-d}\tfrac{(1-qc)(1-qe)(1-qb/a)(qd/a-1)bd}{(1-b)(1-d)(1-qb)(1-qd)qce}
\quad\text{with}\quad
d=\tfrac{qace(1-b)}{qace+qbc+qbe-qbce-q^2bce-ab}.\]
 With this specific value $d$, we can verify that the right member of
 the last equation vanishes. After some simplification, we establish
 the following relation.

\begin{thm}[Two-term contiguous relation of $_3\phi_2$-series]\label{thm-o}
 \bnm
&&_3\phi_2\ffnk{cccc}{q;\frac{qb(1-b)}{qace+qbc+qbe-qbce-q^2bce-ab}}{a,c,e}
{b,\quad\frac{qace(1-b)}{qace+qbc+qbe-qbce-q^2bce-ab}}\\
&=&_3\phi_2\ffnk{cccc}{q;\frac{qb(1-b)}{qace+qbc+qbe-qbce-q^2bce-ab}}{a,q^2c,q^2e}
{q^2b,\frac{q^3ace(1-b)}{qace+qbc+qbe-qbce-q^2bce-ab}}\\
&\times&\tfrac{(qb-a)(1-qc)(1-qe)(qace+qbc+qbe-qbce-q^2ce-ab)}
{(1-qb)(qace+qc+qe-qce-q^2ce-a)(q^2abce+qace+qbc+qbe-q^2ace-q^2bce-qbce-ab)}.
 \enm
\end{thm}

Performing the substitutions $a\to q^a$, $b\to q^b$, $c\to q^c$,
$e\to q^e$ for Theorem \ref{thm-o} and then letting $q\to1$, we
recover the following relation.

\begin{corl}[\citu{wang}{Theorem 15}]
 \bnm
\qquad_3F_2\ffnk{cccc}{1}{a,c,e}{b,\frac{(a-1)(c+e+1)-ce}{b}}
={_3F_2}\ffnk{cccc}{1}{a,c+2,e+2}{b+2,2+\frac{(a-1)(c+e+1)-ce}{b}}\\
\qquad\times\frac{(1+c)(1+e)(a-b-1)\{(c+e-b+1)(a-1)-ce\}}{(b+1)\{(a-1)(c+e+1)-ce\}\{b-ce+(a-1)(c+e+1)\}}.
 \enm
\end{corl}

\subsection{\textbf{C}\&\textbf{D}} $ $\\\\
Let Eq$^{\star}$\eqref{eq-d} stand for Eq\eqref{eq-d} under the
parameter replacements
\[a\to qc,\quad c\to qa,\quad e\to q^2e,\quad b\to q^2b,\quad d\to q^2d.\]
Then consider the linear combination of two equations
\[\text{Eq}\eqref{eq-c}-\text{Eq}^{\star}\eqref{eq-d}\frac{(1-a)(b-c)(c-d)(1-qe)\,bd}{(1-b)(1-d)(1-qb)(1-qd)\,ac^2e}
\quad\text{with}\quad d=\frac{ac(b-1)}{ab+bc-ac-b}.\]
 With this specific value $d$, we can check that the right member of
 the last equation vanishes. After some simplification, we found
 the following relation.

\begin{thm}[Two-term contiguous relation of $_3\phi_2$-series]\label{thm-p}
 \bnm
&&_3\phi_2\ffnk{cccc}{q;\frac{b(b-1)}{(ab+bc-ac-b)e}}{a,c,e}
{b,\quad\frac{ac(b-1)}{ab+bc-ac-b}}\\
&=&_3\phi_2\ffnk{cccc}{q;\frac{b(b-1)}{(ab+bc-ac-b)e}}{qa,qc,q^2e}
{q^2b,\frac{q^2ac(b-1)}{ab+bc-ac-b}}\\
&\times&\frac{(a-b)(b-c)(1-qe)} {e(1-qb)(ab+bc+qac-qabc-ac-b)}.
 \enm
\end{thm}

Employing the substitutions $a\to q^a$, $b\to q^b$, $c\to q^c$,
$e\to q^e$ for Theorem \ref{thm-p} and then letting $q\to1$, we
recover the following relation.

\begin{corl}[\citu{wang}{Theorem 16}]
 \bnm
_3F_2\ffnk{cccc}{1}{a,c,e}{b,\frac{ac}{b}}
={_3F_2}\ffnk{cccc}{1}{a+1,c+1,e+2}{b+2,2+\frac{ac}{b}}
\frac{(a-b)(1+e)(c-b)}{(1+b)(b+ac)}.
 \enm
\end{corl}

\section{Nineteen three-term contiguous relations of $_3\phi_2$-series}
By comparing two patterns of \textbf{A}, \textbf{B}, \textbf{C} and
\textbf{D}, we offer other nineteen three-term contiguous relations
of $_3\phi_2$-series. They produce several two-term contiguous
relations of $_3\phi_2$-series which are different from the ones
before.

\subsection{\textbf{A}\&\textbf{A}}$ $\\\\
Let Eq$^{\star}$\eqref{eq-aa} stand for Eq\eqref{eq-aa} under the
parameter replacements
\[a\to c/q,\quad c\to qa,\quad b\to d/q,\quad d\to qb.\]
Then for an arbitrary variable $Y_q$, the difference
$\text{Eq}\eqref{eq-aa}-Y_q\times\text{Eq}^{\star}\eqref{eq-aa}$
results in the relation:
 \bmn
&&{_3\phi_2}\ffnk{cccc}{q;\frac{bd}{ace}}{a,c,e}{b,d}-Y_q\times
{_3\phi_2}\ffnk{cccc}{q;\frac{bd}{ace}}{qa,c/q,e}{qb,d/q}
 \nnm\\\label{three-aa}
&&=\big(\mathcal{A}_q-\mathcal{A}_q^{\star}Y_q\big)\times
{_4\phi_3}\ffnk{cccc}{q;\frac{bd}{qace}}
{qa,c,e,q\Big(1-\frac{\mathcal{A}_q-\mathcal{A}_q^{\star}Y_q}{\mathfrak{A}_q-\mathfrak{A}_q^{\star}Y_q}\Big)}
{qb,d,\Big(1-\frac{\mathcal{A}_q-\mathcal{A}_q^{\star}Y_q}{\mathfrak{A}_q-\mathfrak{A}_q^{\star}Y_q}\Big)},
 \emn
where the following notations have been used for coefficients
\[\mathcal{A}_q^{\star}=\mathcal{A}_q(c/q,qa,e;\,d/q,qb),\]
\[\mathfrak{A}_q^{\star}=\mathfrak{A}_q(c/q,qa,e;\,d/q,qb)\]
with $\mathcal{A}_q$ and $\mathfrak{A}_q$ being defined in Theorem
\ref{thm-a}. Solving the equation
$1-\frac{\mathcal{A}_q-\mathcal{A}_q^{\star}Y_q}{\mathfrak{A}_q-\mathfrak{A}_q^{\star}Y_q}
=b$ associated with the variable $Y_q$, we obtain from equation
\eqref{three-aa} the following relation.

\begin{thm}[Three-term contiguous relation of $_3\phi_2$-series]\label{thm-q}
 \bnm
_3\phi_2\ffnk{cccc}{q;\frac{bd}{ace}}{a,c,e}{b,d}
=Y_q\times{_3\phi_2}\ffnk{cccc}{q;\frac{bd}{ace}}{qa,c/q,e}{qb,d/q}+
Z_q\times{_3\phi_2}\ffnk{cccc}{q;\frac{bd}{qace}}{qa,c,e}{b,d},
 \enm
where the coefficients $Y_q$ and $Z_q$ are defined by
 \bnm
Y_q&=&
\frac{(b-c)(b-e)(q-d)(qb-c)ad}{(b-1)(qa-d)(q^2abe+acde+b^2d-bde-qabd-qabce)c},\\
Z_q&=&\frac{(qace-bd)(q^2abe+acde+bcd-cde-qabd-qabce)}{(qa-d)(q^2abe+acde+b^2d-bde-qabd-qabce)ce}.
 \enm
\end{thm}

Performing the substitutions $a\to q^a$, $b\to q^b$, $c\to q^c$,
$d\to q^d$, $e\to q^e$ for Theorem \ref{thm-q} and then letting
$q\to1$, we recover the following relation.

\begin{corl}[\citu{wang}{Theorem 19}]
 \bnm
_3F_2\ffnk{cccc}{1}{a,c,e}{b,d}
=Y\,{_3F_2}\ffnk{cccc}{1}{a+1,c-1,e}{b+1,d-1}
+Z\,{_3F_2}\ffnk{cccc}{1}{a+1,c,e}{b,d},
 \enm
where the coefficients $Y$ and $Z$ are given by
 \bnm
Y&=&
\frac{(1+b-c)(b-c)(b-e)(1-d)}{b(1+a-d)(1+b^2+ae+cd+e-c-d-ab-bc-be)},\\
Z&=&\frac{(1+a+c+e-b-d)(1+ae+cd+e-c-d-ab-ce)}{(1+a-d)(1+b^2+ae+cd+e-c-d-ab-bc-be)}.
 \enm
\end{corl}

Specifying the parameter $b\to qa$ in Theorem \ref{thm-q} and using
$q$-Gauss summation formula (cf. \citu{gasper}{p. 14}):
 \bmn
_2\phi_1\ffnk{cccc}{q;\frac{c}{ab}}{a,b}{c}=\ffnk{cccc}{q}{c/a,&c/b}{c,&c/ab}_{\infty}\quad\text{where}\quad
|c/ab|<1, \label{q-gauss}
 \emn
 we get the following relation with
one free parameter less.

\begin{prop}[Two-term contiguous relation of $_3\phi_2$-series]\label{prop-d}
 \bnm
_3\phi_2\ffnk{cccc}{q;\frac{qd}{ce}}{a,c,e}{qa,d}
&=&{_3\phi_2}\ffnk{cccc}{q;\frac{qd}{ce}}{qa,c/q,e}{q^2a,d/q}
\frac{(qa-c)(qa-e)(q-d)(q^2a-c)d}{(qa-1)(q-c)(qa-d)(q^2a-d)ce}\\
&&\xxqdn\qqdn+\:\ffnk{cccc}{q}{d/c,d/e}{d,qd/ce}_{\infty}
 \frac{(q^3a^2e+qacd+acde-cde-q^2a^2d-q^2a^2ce)q}{(q-c)(qa-d)(q^2a-d)e}.
 \enm
\end{prop}

Employing the substitutions $a\to q^a$, $c\to q^c$, $d\to q^d$,
 $e\to q^e$ for Proposition
\ref{prop-d} and then letting $q\to1$, we recover the following
relation.

\begin{corl}[\citu{krattenthaler}{Proposition 1}, see also \citu{wang}{Proposition 20}]
 \bnm
_3F_2\ffnk{cccc}{1}{a,c,e}{a+1,d}&=&{_3F_2}\ffnk{cccc}{1}{a+1,c-1,e}{a+2,d-1}\frac{(1+a-c)(2+a-c)(1+a-e)(d-1)}
{(1+a)(c-1)(1+a-d)(2+a-d)}\\
&+&\Gamma\fnk{cccc}{d,d-c-e+1}{d-c,d-e}
 \frac{(1+ae+cd+e-a-a^2-c-d-ce)}{(c-1)(1+a-d)(2+a-d)},
 \enm
where $\Gamma$-function is given by
 \bnm
  \Gamma(s)=\int_{0}^{\infty}x^{s-1}e^{-x}dx \quad\text{with}\quad
Re(s)>0
 \enm
and the abbreviated expression on $\Gamma$-function is
 \bnm
\Gamma\fnk{cccc}{\alpha,&\beta,&\cdots,&\gamma}{A,&B,&\cdots,&C}
=\frac{\Gamma(\alpha)\Gamma(\beta)\cdots\Gamma(\gamma)}
{\Gamma(A)\Gamma(B)\cdots\Gamma(C)}.
 \enm
\end{corl}
\subsection{\textbf{A}\&\textbf{A}}$ $\\\\
Let Eq$^{\star}$\eqref{eq-a} stand for Eq\eqref{eq-a} under the
parameter replacements
\[a\to c/q,\quad c\to qa.\]
Then for an arbitrary variable $Y_q$, the difference
$\text{Eq}\eqref{eq-a}-Y_q\times\text{Eq}^{\star}\eqref{eq-a}$ leads
us to the relation:
 \bmn
&&_3\phi_2\ffnk{cccc}{q;\frac{bd}{ace}}{a,c,e}{b,d}-Y_q\times
{_3\phi_2}\ffnk{cccc}{q;\frac{bd}{ace}}{qa,c/q,e}{b,d}
  \nnm\\\nnm
&&=\:\big(\mathcal{A}_q-\mathcal{A}_q^{\star}Y_q\big)\times
{_3\phi_2}\ffnk{cccc}{q;\frac{bd}{ace}}{qa,c,e}{qb,d}
\\ \label{three-bbb}
&&+\:\:\big(\mathbb{A}_q-\mathbb{A}_q^{\star}Y_q\big)\times
{_3\phi_2}\ffnk{cccc}{q;\frac{bd}{qace}}{q^2a,qc,qe}{q^2b,qd},
 \emn
where the following notations have been used for coefficients
\[\mathcal{A}_q^{\star}=\mathcal{A}_q(c/q,qa,e;\,b,d),\]
\[\mathbb{A}_q^{\star}=\mathbb{A}_q(c/q,qa,e;\,b,d)\]
with $\mathcal{A}_q$ and $\mathbb{A}_q$ being defined in Theorem
\ref{thm-a}. Solving the equation
$\mathbb{A}_q-\mathbb{A}_q^{\star}Y_q=0$ associated with the
variable $Y_q$, we derive from equation \eqref{three-bbb} the
following relation.

\begin{thm}[Three-term contiguous relation of $_3\phi_2$-series]\label{thm-r}
 \bnm
_3\phi_2\ffnk{cccc}{q;\frac{bd}{ace}}{a,c,e}{b,d}
=Y_q\times{_3\phi_2}\ffnk{cccc}{q;\frac{bd}{ace}}{qa,c/q,e}{b,d}+
Z_q\times{_3\phi_2}\ffnk{cccc}{q;\frac{bd}{ace}}{qa,c,e}{qb,d},
 \enm
where the coefficients $Y_q$ and $Z_q$ are defined by
 \bnm
Y_q=\frac{(c-d)qa}{(qa-d)c},\quad
Z_q=\frac{(b-e)(qa-c)d}{(b-1)(qa-d)ce}.
 \enm
\end{thm}

Performing the substitutions $a\to q^a$, $b\to q^b$, $c\to q^c$,
$d\to q^d$, $e\to q^e$ for Theorem \ref{thm-r} and then letting
$q\to1$, we recover the following relation.

\begin{corl}[\citu{wang}{Theorem 21}]
 \bnm
_3F_2\ffnk{cccc}{1}{a,c,e}{b,d}
=Y\,{_3F_2}\ffnk{cccc}{1}{a+1,c-1,e}{b,d}
+Z\,{_3F_2}\ffnk{cccc}{1}{a+1,c,e}{b+1,d},
 \enm
where the coefficients $Y$ and $Z$ are given by
 \bnm
Y=\frac{c-d}{1+a-d},\quad Z=\frac{(b-e)(1+a-c)}{b(1+a-d)}.
 \enm
\end{corl}

Specifying the parameter $b\to qa$ in Theorem \ref{thm-r} and using
$q$-Gauss summation formula \eqref{q-gauss},
 we deduce the following relation with one free parameter less.

\begin{prop}[Two-term contiguous relation of $_3\phi_2$-series]\label{prop-e}
 \bnm
\quad_3\phi_2\ffnk{cccc}{q;\frac{qd}{ce}}{a,c,e}{qa,d}
={_3\phi_2}\ffnk{cccc}{q;\frac{qd}{ce}}{qa,c,e}{q^2a,d}
\frac{(qa-c)(qa-e)d}{(qa-1)(qa-d)ce}
+\ffnk{cccc}{q}{d/c,d/e}{d,qd/ce}_{\infty}
 \frac{qa}{qa-d}.
 \enm
\end{prop}

Employing the substitutions $a\to q^a$, $c\to q^c$, $d\to q^d$,
$e\to q^e$ for Proposition \ref{prop-e} and then letting $q\to1$, we
recover the following relation.

\begin{corl}[\citu{krattenthaler}{Proposition 3}, see also \citu{wang}{Proposition 22}]
 \bnm
\:\,_3F_2\ffnk{cccc}{1}{a,c,e}{a+1,d}={_3F_2}\ffnk{cccc}{1}{a+1,c,e}{a+2,d}\frac{(1+a-c)(1+a-e)}
{(1+a)(1+a-d)} -\Gamma\fnk{cccc}{d,d-c-e+1}{d-c,d-e}
 \frac{1}{1+a-d}.
 \enm
\end{corl}

Instead, solving the equation
$\mathcal{A}_q-\mathcal{A}_q^{\star}Y_q=0$ associated with the
variable $Y_q$, we attain from equation \eqref{three-bbb} the
following relation.

\begin{thm}[Three-term contiguous relation of $_3\phi_2$-series]\label{thm-s}
 \bnm
_3\phi_2\ffnk{cccc}{q;\frac{bd}{ace}}{a,c,e}{b,d}
=Y_q\times{_3\phi_2}\ffnk{cccc}{q;\frac{bd}{ace}}{qa,c/q,e}{b,d}+
Z_q\times{_3\phi_2}\ffnk{cccc}{q;\frac{bd}{qace}}{q^2a,qc,qe}{q^2b,qd},
 \enm
where the coefficients $Y_q$ and $Z_q$ are defined by
 \bnm
&&Y_q=\frac{(c-d)(qabce+bd+cde-bcd-bde-qace)qa}{(qa-d)(qabce+bd+qade-qabd-bde-qace)c},\\
&&Z_q=\frac{(1-c)(1-e)(1-qa)(c-qa)(b-e)(qace-bd)bd^2}{(1-b)(1-d)(1-qb)(qa-d)(qabce+bd+qade-qabd-bde-qace)qac^2e^2}.
 \enm
\end{thm}

Performing the substitutions $a\to q^a$, $b\to q^b$, $c\to q^c$,
$d\to q^d$, $e\to q^e$ for Theorem \ref{thm-s} and then letting
$q\to1$, we recover the following relation.

\begin{corl}[\citu{wang}{Theorem 23}]
 \bnm
_3F_2\ffnk{cccc}{1}{a,c,e}{b,d}
=Y\,{_3F_2}\ffnk{cccc}{1}{a+1,c-1,e}{b,d}
+Z\,{_3F_2}\ffnk{cccc}{1}{a+2,c+1,e+1}{b+2,d+1},
 \enm
where the coefficients $Y$ and $Z$ are given by
 \bnm
&&Y=\frac{(c-d)(b+ab+ce-bd)}{(1+a-d)(bc+ae+e-bd)},\\
 &&Z=\frac{ce(1+a)(b-e)(1+a-c)(1+a+c+e-b-d)}{bd(1+b)(1+a-d)(bd-bc-ae-e)}.
 \enm
\end{corl}

\textbf{Remark:} There is a tiny mistake in the original equation
due to Chu and Wang \citu{wang}{Theorem 23}. We have added a minus
for the coefficient of the first $_3F_2$-series on the right hand
side.

 Specifying the parameter $b\to qa$ in Theorem \ref{thm-s} and
using $q$-Gauss summation formula \eqref{q-gauss},
 we achieve the following relation with one free parameter less.

\begin{prop}[Two-term contiguous relation of $_3\phi_2$-series]\label{prop-f}
 \bnm
_3\phi_2\ffnk{cccc}{q;\frac{qd}{ce}}{a,c,e}{qa,d}
&=&{_3\phi_2}\ffnk{cccc}{q;\frac{d}{ce}}{q^2a,qc,qe}{q^3a,qd}
\frac{(1-c)(1-e)(qa-c)(qa-e)d^2}{(1-d)(1-qa)(1-q^2a)(qa-d)c^2e^2}\\
&+&\ffnk{cccc}{q}{d/c,d/e}{d,d/ce}_{\infty}
 \frac{(qad+cde+q^2a^2ce-qacd-qace-qade)}{(qa-1)(qa-d)ce}.
 \enm
\end{prop}

Employing the substitutions $a\to q^a$, $c\to q^c$, $d\to q^d$,
$e\to q^e$ for Proposition \ref{prop-f} and then letting $q\to1$, we
recover the following relation.

\begin{corl}[\citu{wang}{Proposition 24}]
 \bnm
_3F_2\ffnk{cccc}{1}{a,c,e}{a+1,d}&=&{_3F_2}\ffnk{cccc}{1}{a+2,c+1,e+1}{a+3,d+1}
\frac{ce(1+a-c)(1+a-e)} {d(1+a)(2+a)(d-a-1)}\\
&+&\Gamma\fnk{cccc}{d,d-c-e}{d-c,d-e}
 \frac{(1+a^2+2a+ce-d-ad)}{(1+a)(1+a-d)}.
 \enm
\end{corl}

Taking $d=\frac{qace(1-qa)}{qa+ce-qac-qae}$ in Proposition
\ref{prop-f}, we establish the following relation.

\begin{prop}[Two-term contiguous relation of $_3\phi_2$-series]\label{prop-ff}
 \bnm
&&_3\phi_2\ffnk{cccc}{q;\frac{q^2a(1-qa)}{qa+ce-qac-qae}}{a,c,e}{qa,\frac{qace(1-qa)}{qa+ce-qac-qae}}\\
&=&{_3\phi_2}\ffnk{cccc}{q;\frac{qa(1-qa)}{qa+ce-qac-qae}}{q^2a,qc,qe}{q^3a,\frac{q^2ace(1-qa)}{qa+ce-qac-qae}}\\
&\times&\frac{(1-qa)(qa-c)(qa-e)}{(1-q^2a)(qa+ce+q^2a^2ce-qac-qae-qace)}.
 \enm
\end{prop}

Performing the substitutions $a\to q^a$, $c\to q^c$, $e\to q^e$
 for Proposition
\ref{prop-ff} and then letting $q\to1$, we recover the following
relation.

\begin{corl}[\citu{wang}{Corollary 25}]
 \bnm
_3F_2\ffnk{cccc}{1}{a,c,e}{a+1,1+a+\frac{ce}{1+a}}&=&{_3F_2}\ffnk{cccc}{1}{a+2,c+1,e+1}{a+3,2+a+\frac{ce}{1+a}}\\
&\times&\frac{(1+a)(1+a-c)(1+a-e)}{(2+a)(1+2a+a^2+ce)}.
 \enm
\end{corl}

\subsection{\textbf{A}\&\textbf{A}}$ $\\\\
Let Eq$^{\star}$\eqref{eq-aa} stand for Eq\eqref{eq-aa} under the
parameter replacements
\[b\to d/q,\quad d\to qb.\]
Then for an arbitrary variable $Y_q$, the difference
$\text{Eq}\eqref{eq-aa}-Y_q\times\text{Eq}^{\star}\eqref{eq-aa}$
results in the relation:
 \bmn
&&{_3\phi_2}\ffnk{cccc}{q;\frac{bd}{ace}}{a,c,e}{b,d}-Y_q\times
{_3\phi_2}\ffnk{cccc}{q;\frac{bd}{ace}}{a,c,e}{qb,d/q}
 \nnm\\\label{three-cc}
&&=\big(\mathcal{A}_q-\mathcal{A}_q^{\star}Y_q\big)\times
{_4\phi_3}\ffnk{cccc}{q;\frac{bd}{qace}}
{qa,c,e,q\Big(1-\frac{\mathcal{A}_q-\mathcal{A}_q^{\star}Y_q}{\mathfrak{A}_q-\mathfrak{A}_q^{\star}Y_q}\Big)}
{qb,d,\Big(1-\frac{\mathcal{A}_q-\mathcal{A}_q^{\star}Y_q}{\mathfrak{A}_q-\mathfrak{A}_q^{\star}Y_q}\Big)},
 \emn
where the following notations have been used for coefficients
\[\mathcal{A}_q^{\star}=\mathcal{A}_q(a,c,e;\,d/q,qb),\]
\[\mathfrak{A}_q^{\star}=\mathfrak{A}_q(a,c,e;\,d/q,qb)\]
with $\mathcal{A}_q$ and $\mathfrak{A}_q$ being defined in Theorem
\ref{thm-a}. Solving the equation
$1-\frac{\mathcal{A}_q-\mathcal{A}_q^{\star}Y_q}{\mathfrak{A}_q-\mathfrak{A}_q^{\star}Y_q}
=b$ associated with the variable $Y_q$, we found from equation
\eqref{three-cc} the following relation.

\begin{thm}[Three-term contiguous relation of $_3\phi_2$-series]\label{thm-t}
 \bnm
_3\phi_2\ffnk{cccc}{q;\frac{bd}{ace}}{a,c,e}{b,d}
=Y_q\times{_3\phi_2}\ffnk{cccc}{q;\frac{bd}{ace}}{a,c,e}{qb,d/q}+
Z_q\times{_3\phi_2}\ffnk{cccc}{q;\frac{bd}{qace}}{qa,c,e}{b,d},
 \enm
where the coefficients $Y_q$ and $Z_q$ are defined by
 \bnm
Y_q&=&
\frac{(a-b)(b-c)(b-e)(d-q)d}{(1-b)(qa-d)(qabce+bcd+bde-b^2d-acde-qbce)},\\
Z_q&=&\frac{(1-a)(d-qb)(qace-bd)}{(qa-d)(qabce+bcd+bde-b^2d-acde-qbce)}.
 \enm
\end{thm}

Employing the substitutions $a\to q^a$, $b\to q^b$, $c\to q^c$,
$d\to q^d$, $e\to q^e$ for Theorem \ref{thm-t} and then letting
$q\to1$, we recover the following relation.

\begin{corl}[\citu{wang}{Theorem 26}]
 \bnm
_3F_2\ffnk{cccc}{1}{a,c,e}{b,d}
=Y\,{_3F_2}\ffnk{cccc}{1}{a,c,e}{b+1,d-1}
+Z\,{_3F_2}\ffnk{cccc}{1}{a+1,c,e}{b,d},
 \enm
where the coefficients $Y$ and $Z$ are given by
 \bnm
Y&=&\frac{(a-b)(b-c)(b-e)(1-d)}{b(1+a-d)(a+ab+bc+be-b^2-ce-ad)},\\
Z&=&\frac{a(1+b-d)(1+a+c+e-b-d)}{(1+a-d)(a+ab+bc+be-b^2-ce-ad)}.
 \enm
\end{corl}

Specifying the parameter $b\to qa$ in Theorem \ref{thm-t} and using
$q$-Gauss summation formula \eqref{q-gauss},
 we obtain the following relation with
one free parameter less.

\begin{prop}[Two-term contiguous relation of $_3\phi_2$-series] \label{prop-g}
 \bnm
_3\phi_2\ffnk{cccc}{q;\frac{qd}{ce}}{a,c,e}{qa,d}
&=&{_3\phi_2}\ffnk{cccc}{q;\frac{qd}{ce}}{a,c,e}{q^2a,d/q}
\frac{(1-q)(q-d)(qa-c)(qa-e)d}{\sst(1-qa)(qa-d)(q^2ad+q^2ce+cde-qcd-qde-q^2ace)}\\
&+&\ffnk{cccc}{q}{d/c,d/e}{d,d/ce}_{\infty}
\frac{(1-a)(ce-d)(q^2a-d)q}{(qa-d)(q^2ad+q^2ce+cde-qcd-qde-q^2ace)}.
 \enm
\end{prop}

Performing the substitutions $a\to q^a$, $c\to q^c$, $d\to q^d$,
$e\to q^e$ for Proposition \ref{prop-g} and then letting $q\to1$, we
recover the following relation.

\begin{corl}[\citu{wang}{Proposition 27}]
 \bnm
_3F_2\ffnk{cccc}{1}{a,c,e}{a+1,d}&=&{_3F_2}\ffnk{cccc}{1}{a,c,e}{a+2,d-1}
\frac{(1+a-c)(1+a-e)(d-1)}
{(1+a)(1+a-d)(ac+ae+c+e-ce-ad-1)}\\
&+&\Gamma\fnk{cccc}{d,d-c-e}{d-c,d-e}
 \frac{a(2+a-d)(c+e-d)}{(1+a-d)(ac+ae+c+e-ce-ad-1)}.
 \enm
\end{corl}
\subsection{\textbf{A}\&\textbf{A}}
Let Eq$^{\star}$\eqref{eq-a} stand for Eq\eqref{eq-a} under the
parameter replacements
\[b\to d/q,\quad d\to qb.\]
Then for an arbitrary variable $Y_q$, the difference
$\text{Eq}\eqref{eq-a}-Y_q\times\text{Eq}^{\star}\eqref{eq-a}$ leads
us to the relation:
 \bmn
&&_3\phi_2\ffnk{cccc}{q;\frac{bd}{ace}}{a,c,e}{b,d}-Y_q\times
{_3\phi_2}\ffnk{cccc}{q;\frac{bd}{ace}}{a,c,e}{qb,d/q}
  \nnm\\\nnm
&&=\:\big(\mathcal{A}_q-\mathcal{A}_q^{\star}Y_q\big)\times
{_3\phi_2}\ffnk{cccc}{q;\frac{bd}{ace}}{qa,c,e}{qb,d}
\\ \label{three-ddd}
&&+\:\:\big(\mathbb{A}_q-\mathbb{A}_q^{\star}Y_q\big)\times
{_3\phi_2}\ffnk{cccc}{q;\frac{bd}{qace}}{q^2a,qc,qe}{q^2b,qd},
 \emn
where the following notations have been used for coefficients
\[\mathcal{A}_q^{\star}=\mathcal{A}_q(a,c,e;\,d/q,qb),\]
\[\mathbb{A}_q^{\star}=\mathbb{A}_q(a,c,e;\,d/q,qb)\]
with $\mathcal{A}_q$ and $\mathbb{A}_q$ being defined in Theorem
\ref{thm-a}. Solving the equation
$\mathbb{A}_q-\mathbb{A}_q^{\star}Y_q=0$ associated with the
variable $Y_q$, we get from equation \eqref{three-ddd} the following
relation.

\begin{thm}[Three-term contiguous relation of $_3\phi_2$-series]\label{thm-u}
 \bnm
_3\phi_2\ffnk{cccc}{q;\frac{bd}{ace}}{a,c,e}{b,d}
=Y_q\times{_3\phi_2}\ffnk{cccc}{q;\frac{bd}{ace}}{a,c,e}{qb,d/q}+
Z_q\times{_3\phi_2}\ffnk{cccc}{q;\frac{bd}{ace}}{qa,c,e}{qb,d},
 \enm
where the coefficients $Y_q$ and $Z_q$ are defined by
 \bnm
Y_q=\frac{(a-b)(q-d)}{(1-b)(qa-d)},\quad
Z_q=\frac{(1-a)(qb-d)}{(1-b)(qa-d)}.
 \enm
\end{thm}

Employing the substitutions $a\to q^a$, $b\to q^b$, $c\to q^c$,
$d\to q^d$, $e\to q^e$ for Theorem \ref{thm-u} and then letting
$q\to1$, we recover the following relation.

\begin{corl}[\citu{wang}{Theorem 28}]
 \bnm
_3F_2\ffnk{cccc}{1}{a,c,e}{b,d}
=Y\,{_3F_2}\ffnk{cccc}{1}{a,c,e}{b+1,d-1}
+Z\,{_3F_2}\ffnk{cccc}{1}{a+1,c,e}{b+1,d},
 \enm
where the coefficients $Y$ and $Z$ are given by
 \bnm
Y=\frac{(a-b)(d-1)}{b(1+a-d)},\quad Z=\frac{a(1+b-d)}{b(1+a-d)}.
 \enm
\end{corl}
Specifying the parameter $d\to qc$ in Theorem \ref{thm-u} and using
$q$-Gauss summation formula \eqref{q-gauss},
 we derive the following relation with one free parameter less under the replacements $a\to c, c\to a$.

\begin{prop} [Two-term contiguous relation of $_3\phi_2$-series]\label{prop-h}
 \bnm
\quad_3\phi_2\ffnk{cccc}{q;\frac{qb}{ce}}{a,c,e}{qa,b}
={_3\phi_2}\ffnk{cccc}{q;\frac{qb}{ce}}{a,qc,e}{qa,qb}
\frac{(a-b)(1-c)}{(a-c)(1-b)}
+\ffnk{cccc}{q}{b/c,qb/e}{b,qb/ce}_{\infty}
 \frac{(a-1)c}{a-c}.
 \enm
\end{prop}

Performing the substitutions $a\to q^a$,  $b\to q^b$, $c\to q^c$,
$e\to q^e$ for Proposition \ref{prop-h} and then letting $q\to1$, we
recover the following relation.

\begin{corl}[\citu{wang}{Proposition 29}]
 \bnm
\:\,_3F_2\ffnk{cccc}{1}{a,c,e}{a+1,b}={_3F_2}\ffnk{cccc}{1}{a,c+1,e}{a+1,b+1}\frac{(a-b)c}
{(a-c)b}+\Gamma\fnk{cccc}{b,b-c-e+1}{b-c,b-e+1}
 \frac{a}{a-c}.
 \enm
\end{corl}

Instead, solving the equation
$\mathcal{A}_q-\mathcal{A}_q^{\star}Y_q=0$ associated with the
variable $Y_q$, we deduce from equation \eqref{three-ddd} the
following relation.

\begin{thm}[Three-term contiguous relation of $_3\phi_2$-series]\label{thm-v}
 \bnm
_3\phi_2\ffnk{cccc}{q;\frac{bd}{ace}}{a,c,e}{b,d}
=Y_q\times{_3\phi_2}\ffnk{cccc}{q;\frac{bd}{ace}}{a,c,e}{qb,d/q}+
Z_q\times{_3\phi_2}\ffnk{cccc}{q;\frac{bd}{qace}}{q^2a,qc,qe}{q^2b,qd},
 \enm
where the coefficients $Y_q$ and $Z_q$ are defined by
 \bnm
&&Y_q=\frac{(a-b)(q-d)(qabce+cde+bd-bcd-bde-qace)}{(1-b)(qa-d)(qbce+acde+bd-bcd-bde-qace)},\\
&&Z_q=\frac{(1-a)(1-c)(1-e)(1-qa)(qb-d)(qace-bd)bd}{(1-b)(1-d)(1-qb)(qa-d)(qbce+acde+bd-bcd-bde-qace)qace}.
 \enm
\end{thm}

Employing the substitutions $a\to q^a$, $b\to q^b$, $c\to q^c$,
$d\to q^d$, $e\to q^e$ for Theorem \ref{thm-v} and then letting
$q\to1$, we recover the following relation.

\begin{corl}[\citu{wang}{Theorem 30}]
 \bnm
_3F_2\ffnk{cccc}{1}{a,c,e}{b,d}
=Y\,{_3F_2}\ffnk{cccc}{1}{a,c,e}{b+1,d-1}
+Z\,{_3F_2}\ffnk{cccc}{1}{a+2,c+1,e+1}{b+2,d+1},
 \enm
where the coefficients $Y$ and $Z$ are given by
 \bnm
&&Y=\frac{(a-b)(d-1)(bd-ab-ce-b)}{b(1+a-d)(a+bd-ad-ce-b)},\\
 &&Z=\frac{ace(1+a)(1+b-d)(1+a+c+e-b-d)}{bd(1+b)(1+a-d)(a+bd-ad-ce-b)}.
 \enm
\end{corl}

 Specifying the parameter $b\to c$ in Theorem \ref{thm-v} and
using $q$-Gauss summation formula \eqref{q-gauss},
 we attain the following relation with one free parameter less under the replacements
 $a\to c,\: c\to a, \:d\to qd$.

\begin{prop}[Two-term contiguous relation of $_3\phi_2$-series]\label{prop-i}
 \bnm
_3\phi_2\ffnk{cccc}{q;\frac{qd}{ce}}{a,c,e}{qa,d}
&=&{_3\phi_2}\ffnk{cccc}{q;\frac{d}{ce}}{qa,q^2c,qe}{q^2a,q^2d}
\frac{(1-c)(1-e)(1-qc)(a-d)d}{(1-d)(1-qd)(1-qa)(c-a)ce}\\
&+&\ffnk{cccc}{q}{d/c,qd/e}{d,d/ce}_{\infty}
 \frac{(ae+cde+d-ad-de-ce)}{(a-c)e}.
 \enm
\end{prop}

Performing the substitutions $a\to q^a$, $c\to q^c$, $d\to q^d$,
$e\to q^e$ for Proposition \ref{prop-i} and then letting $q\to1$, we
recover the following relation.

\begin{corl}[\citu{wang}{Proposition 31}]
 \bnm
_3F_2\ffnk{cccc}{1}{a,c,e}{a+1,d}&=&{_3F_2}\ffnk{cccc}{1}{a+1,c+2,e+1}{a+2,d+2}
\frac{ce(1+c)(a-d)}{d(1+a)(1+d)(c-a)}\\
&+&\Gamma\fnk{cccc}{d,d-c-e}{d-c,d-e+1}
 \frac{(ad-cd-ae)}{(a-c)}.
 \enm
\end{corl}

Taking $d=\frac{e(c-a)}{1+ce-a-e}$ in Proposition \ref{prop-i}, we
achieve the following relation.

\begin{prop}[Two-term contiguous relation of $_3\phi_2$-series]\label{prop-j}
 \bnm
&&_3\phi_2\ffnk{cccc}{q;\frac{q(c-a)}{c(1+ce-a-e)}}{a,c,e}{qa,\frac{e(c-a)}{1+ce-a-e}}\\
&=&_3\phi_2\ffnk{cccc}{q;\frac{c-a}{c(1+ce-a-e)}}{qa,q^2c,qe}{q^2a,\frac{q^2e(c-a)}{1+ce-a-e}}\\
&\times&\frac{(1-c)(1-qc)(a-ce)}{(1-qa)(1+ce+qae-qce-a-e)c}.
 \enm
\end{prop}

Employing the substitutions $a\to q^a$, $c\to q^c$, $e\to q^e$
 for Proposition
\ref{prop-j} and then letting $q\to1$, we recover the following
relation.

\begin{corl}[\citu{wang}{Corollary 32}]
 \bnm
_3F_2\ffnk{cccc}{1}{a,c,e}{a+1,\frac{ae}{a-c}}={_3F_2}\ffnk{cccc}{1}{a+1,c+2,e+1}{a+2,2+\frac{ae}{a-c}}
\frac{c(1+c)(c+e-a)}{(1+a)(a+ae-c)}.
 \enm
\end{corl}
\subsection{\textbf{A}\&\textbf{B}}$ $\\\\
Let Eq$^{\star}$\eqref{eq-bb} stand for Eq\eqref{eq-bb} under the
parameter replacements
\[a\to qc,\quad c\to qa,\quad b\to d,\quad d\to q^2b.\]
Then for an arbitrary variable $Y_q$, the difference
$\text{Eq}\eqref{eq-aa}-Y_q\times\text{Eq}^{\star}\eqref{eq-bb}$
results in the relation:
 \bmn
&&{_3\phi_2}\ffnk{cccc}{q;\frac{bd}{ace}}{a,c,e}{b,d}-Y_q\times
{_3\phi_2}\ffnk{cccc}{q;\frac{bd}{ace}}{qa,qc,e}{q^2b,d}
  \nnm\\\label{three-ee}
&&=\big(\mathcal{A}_q-\mathcal{B}_q^{\star}Y_q\big)\times
{_4\phi_3}\ffnk{cccc}{q;\frac{bd}{qace}}
{qa,c,e,q\Big(1-\frac{\mathcal{A}_q-\mathcal{B}_q^{\star}Y_q}{\mathfrak{A}_q-\mathfrak{B}_q^{\star}Y_q}\Big)}
{qb,d,\Big(1-\frac{\mathcal{A}_q-\mathcal{B}_q^{\star}Y_q}{\mathfrak{A}_q-\mathfrak{B}_q^{\star}Y_q}\Big)},
 \emn
where the following notations have been used for coefficients
\[\mathcal{B}_q^{\star}=\mathcal{B}_q(qc,qa,e;\,d,q^2b),\]
\[\mathfrak{B}_q^{\star}=\mathfrak{B}_q(qc,qa,e;\,d,q^2b)\]
with $\mathcal{B}_q$ and $\mathfrak{B}_q$ being defined in Theorem
\ref{thm-b}. Solving the equation
$1-\frac{\mathcal{A}_q-\mathcal{B}_q^{\star}Y_q}{\mathfrak{A}_q-\mathfrak{B}_q^{\star}Y_q}
=e$ associated with the variable $Y_q$, we establish from equation
\eqref{three-ee} the following relation.

\begin{thm}[Three-term contiguous relation of $_3\phi_2$-series]\label{thm-w}
 \bnm
_3\phi_2\ffnk{cccc}{q;\frac{bd}{ace}}{a,c,e}{b,d}
=Y_q\times{_3\phi_2}\ffnk{cccc}{q;\frac{bd}{ace}}{qa,qc,e}{q^2b,d}+
Z_q\times{_3\phi_2}\ffnk{cccc}{q;\frac{bd}{qace}}{qa,c,qe}{qb,d},
 \enm
where the coefficients $Y_q$ and $Z_q$ are defined by
 \bnm
Y_q&=&
\frac{(1-c)(b-a)(b-e)(qb-e)(qae-d)d}{(1-b)(1-qb)(qa-d)(qace+acde+bde-ade-bcd-qace^2)e},\\
Z_q&=&\frac{(1-e)(qace-bd)(qabc+ad+cd-bd-acd-qac)}{(1-b)(d-qa)(qace+acde+bde-ade-bcd-qace^2)c}.
 \enm
\end{thm}

Performing the substitutions $a\to q^a$, $b\to q^b$, $c\to q^c$,
$d\to q^d$, $e\to q^e$ for Theorem \ref{thm-w} and then letting
$q\to1$, we recover the following relation.

\begin{corl}[\citu{wang}{Theorem 33}]
 \bnm
_3F_2\ffnk{cccc}{1}{a,c,e}{b,d}
=Y\,{_3F_2}\ffnk{cccc}{1}{a+1,c+1,e}{b+2,d}
+Z\,{_3F_2}\ffnk{cccc}{1}{a+1,c,e+1}{b+1,d},
 \enm
where the coefficients $Y$ and $Z$ are given by
 \bnm
Y&=&
\frac{c(a-b)(b-e)(1+b-e)(1+a+e-d)}{b(1+b)(1+a-d)(ac+be+de-ae-e-e^2-bc)},\\
Z&=&\frac{e(b+d-a-c-e-1)(b+ab+bc-bd-ac)}{b(1+a-d)(ac+be+de-ae-e-e^2-bc)}.
 \enm
\end{corl}

Specifying the parameter $d\to qe$ in Theorem \ref{thm-w} and using
$q$-Gauss summation formula \eqref{q-gauss},
 we found the following relation with
one free parameter less.

\begin{prop}[Two-term contiguous relation of $_3\phi_2$-series]\label{prop-k}
 \bnm
_3\phi_2\ffnk{cccc}{q;\frac{qb}{ac}}{a,c,e}{b,qe}
&=&{_3\phi_2}\ffnk{cccc}{q;\frac{qb}{ac}}{qa,qc,e}{q^2b,qe}
\frac{(1-a)(1-c)(b-e)(qb-e)}{(1-b)(1-qb)(a-e)(c-e)}\\
&&\xxqdn\qqdn+\:\ffnk{cccc}{q}{qb/a,qb/c}{b,qb/ac}_{\infty}
 \frac{(e-1)(abc+ae+ce-ac-be-ace)}{(a-e)(c-e)}.
 \enm
\end{prop}

Employing the substitutions $a\to q^a$, $b\to q^b$, $c\to q^c$,
$e\to q^e$ for Proposition \ref{prop-k} and then letting $q\to1$, we
recover the following relation.

\begin{corl}[\citu{krattenthaler}{Proposition 2}, see also \citu{wang}{Proposition 34}]
 \bnm
_3F_2\ffnk{cccc}{1}{a,c,e}{b,e+1}&=&{_3F_2}\ffnk{cccc}{1}{a+1,c+1,e}{b+2,e+1}\frac{ac(b-e)(1+b-e)}
{b(1+b)(a-e)(c-e)}\\
&+&\Gamma\fnk{cccc}{b,b-a-c+1}{b-a+1,b-c+1}
 \frac{e(ac+be-ab-bc)}{(a-e)(c-e)}.
 \enm
\end{corl}

Instead, solving the equation
$1-\frac{\mathcal{A}_q-\mathcal{B}_q^{\star}Y_q}{\mathfrak{A}_q-\mathfrak{B}_q^{\star}Y_q}
=b$ associated with the variable $Y_q$, we obtain from equation
\eqref{three-ee} the following relation.

\begin{thm}[Three-term contiguous relation of $_3\phi_2$-series]\label{thm-x}
 \bnm
_3\phi_2\ffnk{cccc}{q;\frac{bd}{ace}}{a,c,e}{b,d}
=Y_q\times{_3\phi_2}\ffnk{cccc}{q;\frac{bd}{ace}}{qa,qc,e}{q^2b,d}+
Z_q\times{_3\phi_2}\ffnk{cccc}{q;\frac{bd}{qace}}{qa,c,e}{b,d},
 \enm
where the coefficients $Y_q$ and $Z_q$ are defined by
 \bnm
Y_q&=&
\frac{(a-b)(b-c)(1-c)(b-e)(qb-e)d^2}{(1-b)(1-qb)(qa-d)(qabce+ade+bcd-b^2d-acde-qace)ce},\\
Z_q&=&\frac{(qace-bd)(qabc+ad+cd-bd-acd-qac)}{(qa-d)(qabce+ade+bcd-b^2d-acde-qace)c}.
 \enm
\end{thm}

Performing the substitutions $a\to q^a$, $b\to q^b$, $c\to q^c$,
$d\to q^d$, $e\to q^e$ for Theorem \ref{thm-x} and then letting
$q\to1$, we recover the following relation.

\begin{corl}[\citu{wang}{Theorem 35}]
 \bnm
_3F_2\ffnk{cccc}{1}{a,c,e}{b,d}
=Y\,{_3F_2}\ffnk{cccc}{1}{a+1,c+1,e}{b+2,d}
+Z\,{_3F_2}\ffnk{cccc}{1}{a+1,c,e}{b,d},
 \enm
where the coefficients $Y$ and $Z$ are given by
 \bnm
Y&=&
\frac{(a-b)(b-c)(b-e)(1+b-e)c}{b(1+b)(1+a-d)(ac+ce+bd+b^2-b-ab-be-2bc)},\\
Z&=&\frac{(b+d-a-c-e-1)(ab+bc+b-bd-ac)}{(1+a-d)(ac+ce+bd+b^2-b-ab-be-2bc)}.
 \enm
\end{corl}

\textbf{Remark:} There is a tiny mistake in the original equation
due to Chu and Wang \citu{wang}{Theorem 35}. We have changed the
parameter $e+1$ into $e$ for the last $_3F_2$-series. The similar
changes should also be used in Chu and Wang \citu{wang}{Proposition
36 and Corollary 37}.

Specifying the parameter $b\to qa$ in Theorem \ref{thm-x} and using
$q$-Gauss summation formula \eqref{q-gauss},
 we get the following relation with
one free parameter less.

\begin{prop}[Two-term contiguous relation of $_3\phi_2$-series]\label{prop-l}
 \bnm
&&\!\!\xqdn_3\phi_2\ffnk{cccc}{q;\frac{qd}{ce}}{a,c,e}{qa,d}
={_3\phi_2}\ffnk{cccc}{q;\frac{qd}{ce}}{qa,qc,e}{q^3a,d}
\tfrac{(1-q)(1-c)(qa-c)(qa-e)(q^2a-e)d^2}{(qa-d)(qa-1)(q^2a-1)(q^2ace+qcd+de-cde-qce-q^2ad)ce}\\
&&\qquad\qquad\qquad\quad\qdn+\:\:\ffnk{cccc}{q}{d/c,d/e}{d,qd/ce}_{\infty}
\frac{qe(q^2a^2c+ad+cd-acd-qac-qad)}{(qa-d)(q^2ace+qcd+de-cde-qce-q^2ad)}.
 \enm
\end{prop}

Employing the substitutions $a\to q^a$, $c\to q^c$, $d\to q^d$,
$e\to q^e$ for Proposition \ref{prop-l} and then letting $q\to1$, we
recover the following relation.

\begin{corl}[\citu{wang}{Proposition 36}]
 \bnm
_3F_2\ffnk{cccc}{1}{a,c,e}{a+1,d}&=&{_3F_2}\ffnk{cccc}{1}{a+1,c+1,e}{a+3,d}
 \tfrac{c(1+a-c)(1+a-e)(2+a-e)}{(a+1)(a+2)(1+a-d)(ac+ae+2c+e-d-ad-ce)}\\
&+&\Gamma\fnk{cccc}{d,d-c-e+1}{d-c,d-e}
 \frac{(1+2a+a^2+c-d-ad)}{(d-a-1)(ac+ae+2c+e-d-ad-ce)}.
 \enm
\end{corl}

Taking $d=\frac{qac(qa-1)}{qa+ac-a-c}$ in Proposition \ref{prop-l},
we derive the following relation.

\begin{prop}[Two-term contiguous relation of $_3\phi_2$-series]\label{prop-m}
 \bnm
&&\xqdn_3\phi_2\ffnk{cccc}{q;\frac{q^2a(qa-1)}{e(qa+ac-a-c)}}{a,c,e}{qa,\:\frac{qac(qa-1)}{qa+ac-a-c}}
={_3\phi_2}\ffnk{cccc}{q;\frac{q^2a(qa-1)}{e(qa+ac-a-c)}}{qa,qc,e}{q^3a,\:\frac{qac(qa-1)}{qa+ac-a-c}}
\\&&\qquad\qquad\qquad\qquad\qquad\qquad\qquad\qquad\qquad\!\!\times\:\frac{(q^2a-e)}{(q^2a-1)e}.
 \enm
\end{prop}

Performing the substitutions $a\to q^a$, $c\to q^c$, $e\to q^e$ for
 Proposition \ref{prop-m} and then
letting $q\to1$, we recover the following relation.

\begin{corl}[\citu{wang}{Corollary 37}]
 \bnm
_3F_2\ffnk{cccc}{1}{a,c,e}{a+1,1+a+\frac{c}{1+a}}&=&{_3F_2}\ffnk{cccc}{1}{a+1,c+1,e}{a+3,1+a+\frac{c}{1+a}}
 \frac{(2+a-e)}{(2+a)}.
 \enm
\end{corl}
\subsection{\textbf{A}\&\textbf{B}}$ $\\\\
Let Eq$^{\star}$\eqref{eq-bb} stand for Eq\eqref{eq-bb} under the
parameter replacements
\[a\to qc,\quad c\to qa,\quad b\to qb,\quad d\to qd.\]
Then for an arbitrary variable $Y_q$, the difference
$\text{Eq}\eqref{eq-aa}-Y_q\times\text{Eq}^{\star}\eqref{eq-bb}$
leads us to the relation:
 \bmn
&&{_3\phi_2}\ffnk{cccc}{q;\frac{bd}{ace}}{a,c,e}{b,d}-Y_q\times
{_3\phi_2}\ffnk{cccc}{q;\frac{bd}{ace}}{qa,qc,e}{qb,qd}
  \nnm\\\label{three-ff}
&&=\big(\mathcal{A}_q-\mathcal{B}_q^{\star}Y_q\big)\times
{_4\phi_3}\ffnk{cccc}{q;\frac{bd}{qace}}
{qa,c,e,q\Big(1-\frac{\mathcal{A}_q-\mathcal{B}_q^{\star}Y_q}{\mathfrak{A}_q-\mathfrak{B}_q^{\star}Y_q}\Big)}
{qb,d,\Big(1-\frac{\mathcal{A}_q-\mathcal{B}_q^{\star}Y_q}{\mathfrak{A}_q-\mathfrak{B}_q^{\star}Y_q}\Big)},
 \emn
where the following notations have been used for coefficients
\[\mathcal{B}_q^{\star}=\mathcal{B}_q(qc,qa,e;\,qb,qd),\]
\[\mathfrak{B}_q^{\star}=\mathfrak{B}_q(qc,qa,e;\,qb,qd)\]
with $\mathcal{B}_q$ and $\mathfrak{B}_q$ being defined in Theorem
\ref{thm-b}. Solving the equation
$1-\frac{\mathcal{A}_q-\mathcal{B}_q^{\star}Y_q}{\mathfrak{A}_q-\mathfrak{B}_q^{\star}Y_q}
=b$ associated with the variable $Y_q$, we deduce from equation
\eqref{three-ff} the following relation.

\begin{thm}[Three-term contiguous relation of $_3\phi_2$-series]\label{thm-y}
 \bnm
_3\phi_2\ffnk{cccc}{q;\frac{bd}{ace}}{a,c,e}{b,d}
=Y_q\times{_3\phi_2}\ffnk{cccc}{q;\frac{bd}{ace}}{qa,qc,e}{qb,qd}+
Z_q\times{_3\phi_2}\ffnk{cccc}{q;\frac{bd}{qace}}{qa,c,e}{b,d},
 \enm
where the coefficients $Y_q$ and $Z_q$ are defined by
 \bnm
Y_q= \frac{(1-c)(b-e)(d-e)bd}{(1-b)(1-d)(qae-bd)ce},\quad
Z_q=\frac{(qace-bd)}{(qae-bd)c}.
 \enm
\end{thm}

Employing the substitutions $a\to q^a$, $b\to q^b$, $c\to q^c$,
$d\to q^d$, $e\to q^e$ for Theorem \ref{thm-y} and then letting
$q\to1$, we recover the following relation.

\begin{corl}[\citu{wang}{Theorem 38}]
 \bnm
_3F_2\ffnk{cccc}{1}{a,c,e}{b,d}
=Y\,{_3F_2}\ffnk{cccc}{1}{a+1,c+1,e}{b+1,d+1}
+Z\,{_3F_2}\ffnk{cccc}{1}{a+1,c,e}{b,d},
 \enm
where the coefficients $Y$ and $Z$ are given by
 \bnm
Y=\frac{c(b-e)(e-d)}{bd(1+a+e-b-d)},\quad
Z=\frac{(1+a+c+e-b-d)}{(1+a+e-b-d)}.
 \enm
\end{corl}

Specifying the parameter $b\to qa$ in Theorem \ref{thm-y} and using
$q$-Gauss summation formula \eqref{q-gauss},
 we attain the following relation with
one free parameter less.

\begin{prop}[Two-term contiguous relation of $_3\phi_2$-series]\label{prop-n}
 \bnm
_3\phi_2\ffnk{cccc}{q;\frac{qd}{ce}}{a,c,e}{qa,d}
={_3\phi_2}\ffnk{cccc}{q;\frac{qd}{ce}}{qa,qc,e}{q^2a,qd}
\frac{(qa-e)(1-c)d}{(qa-1)(1-d)ce}+\ffnk{cccc}{q}{d/c,qd/e}{d,qd/ce}_{\infty}.
 \enm
\end{prop}

Performing the substitutions $a\to q^a$, $c\to q^c$, $d\to q^d$,
$e\to q^e$ for Proposition \ref{prop-n} and then letting $q\to1$, we
recover the following relation.

\begin{corl}[\citu{wang}{Proposition 39}]
 \bnm
_3F_2\ffnk{cccc}{1}{a,c,e}{a+1,d}={_3F_2}\ffnk{cccc}{1}{a+1,c+1,e}{a+2,d+1}
 \frac{(1+a-e)c}{(1+a)d}+\Gamma\fnk{cccc}{d,d-c-e+1}{d-c,d-e+1}.
 \enm
\end{corl}

\subsection{\textbf{A}\&\textbf{D}}$ $\\\\
Let Eq$^{\star}$\eqref{eq-d} stand for Eq\eqref{eq-d} under the
parameter replacements
\[a\to qa,\quad c\to qc,\quad e\to qe,\quad b\to q^2b,\quad d\to qd.\]
Then for an arbitrary variable $Y_q$, the difference
$\text{Eq}\eqref{eq-a}-Y_q\times\text{Eq}^{\star}\eqref{eq-d}$
results in relation:
 \bmn
&&_3\phi_2\ffnk{cccc}{q;\frac{bd}{ace}}{a,c,e}{b,d}-Y_q\times
{_3\phi_2}\ffnk{cccc}{q;\frac{bd}{ace}}{qa,qc,qe}{q^2b,qd}
  \nnm\\\nnm
&&=\:\big(\mathcal{A}_q-\mathcal{D}_q^{\star}Y_q\big)\times
{_3\phi_2}\ffnk{cccc}{q;\frac{bd}{ace}}{qa,c,e}{qb,d}
\\ \label{three-ggg}
&&+\:\:\big(\mathbb{A}_q-\mathbb{D}_q^{\star}Y_q\big)\times
{_3\phi_2}\ffnk{cccc}{q;\frac{bd}{qace}}{q^2a,qc,qe}{q^2b,qd},
 \emn
where the following notations have been used for coefficients
\[\mathcal{D}_q^{\star}=\mathcal{D}_q(qa,qc,qe;\,q^2b,qd),\]
\[\mathbb{D}_q^{\star}=\mathbb{D}_q(qa,qc,qe;\,q^2b,qd)\]
with $\mathcal{D}_q$ and $\mathbb{D}_q$ being defined in Theorem
\ref{thm-d}. Solving the equation
$\mathbb{A}_q-\mathbb{D}_q^{\star}Y_q=0$ associated with the
variable $Y_q$, we achieve from equation \eqref{three-ggg} the
following relation.

\begin{thm}[Three-term contiguous relation of $_3\phi_2$-series]\label{thm-z}
 \bnm
\:\:_3\phi_2\ffnk{cccc}{q;\frac{bd}{ace}}{a,c,e}{b,d}
={_3\phi_2}\!\ffnk{cccc}{q;\frac{bd}{ace}}{qa,qc,qe}{q^2b,qd}\frac{(a-b)(1-c)(1-e)bd}{(1-b)(1-qb)(d-1)ace}+
{_3\phi_2}\!\ffnk{cccc}{q;\frac{bd}{ace}}{qa,c,e}{qb,d}.
 \enm
\end{thm}

Employing the substitutions $a\to q^a$, $b\to q^b$, $c\to q^c$,
$d\to q^d$, $e\to q^e$ for Theorem \ref{thm-z} and then letting
$q\to1$, we recover the following relation.

\begin{corl}[\citu{wang}{Theorem 40}]
 \bnm
_3F_2\ffnk{cccc}{1}{a,c,e}{b,d}
={_3F_2}\ffnk{cccc}{1}{a+1,c+1,e+1}{b+2,d+1}\frac{(a-b)ce}{(1+b)bd}
+{_3F_2}\ffnk{cccc}{1}{a+1,c,e}{b+1,d}.
 \enm
\end{corl}

Specifying the parameter $d\to qa$ in Theorem \ref{thm-z} and using
$q$-Gauss summation formula \eqref{q-gauss},
 we establish the following relation with one free parameter less.

\begin{prop} [Two-term contiguous relation of $_3\phi_2$-series]\label{prop-o}
 \bnm
\quad_3\phi_2\ffnk{cccc}{q;\frac{qb}{ce}}{a,c,e}{qa,b}
={_3\phi_2}\ffnk{cccc}{q;\frac{qb}{ce}}{qa,qc,qe}{q^2a,q^2b}
\frac{(a-b)(1-c)(1-e)qb}{(1-b)(1-qb)(qa-1)ce}
+\ffnk{cccc}{q}{qb/c,qb/e}{qb,qb/ce}_{\infty}.
 \enm
\end{prop}

Performing the substitutions $a\to q^a$, $b\to q^b$, $c\to q^c$,
$e\to q^e$ for Proposition \ref{prop-o} and then letting $q\to1$, we
recover the following relation.

\begin{corl}[\citu{wang}{Proposition 41}]
 \bnm
\qquad_3F_2\ffnk{cccc}{1}{a,c,e}{a+1,b}={_3F_2}\ffnk{cccc}{1}{a+1,c+1,e+1}{a+2,b+2}\frac{(a-b)ce}{(1+a)(1+b)b}
+\Gamma\fnk{cccc}{b+1,b-c-e+1}{b-c+1,b-e+1}.
 \enm
\end{corl}

\subsection{\textbf{B}\&\textbf{B}}$ $\\\\
Let Eq$^{\star}$\eqref{eq-b} stand for Eq\eqref{eq-b} under the
parameter replacements
\[a\to qc, \quad c\to a/q.\]
Then for an arbitrary variable $Y_q$, the difference
$\text{Eq}\eqref{eq-b}-Y_q\times\text{Eq}^{\star}\eqref{eq-b}$ leads
us to the relation:
 \bmn
&&_3\phi_2\ffnk{cccc}{q;\frac{bd}{ace}}{a,c,e}{b,d}-Y_q\times
{_3\phi_2}\ffnk{cccc}{q;\frac{bd}{ace}}{a/q,qc,e}{b,d}
  \nnm\\\nnm
&&=\:\big(\mathcal{B}_q-\mathcal{B}_q^{\star}Y_q\big)\times
{_3\phi_2}\ffnk{cccc}{q;\frac{bd}{ace}}{a/q,c,e}{b,d/q}
\\ \label{three-hhh}
&&+\:\:\big(\mathbb{B}_q-\mathbb{B}_q^{\star}Y_q\big)\times
{_3\phi_2}\ffnk{cccc}{q;\frac{bd}{qace}}{a,qc,qe}{qb,d},
 \emn
where the following notations have been used for coefficients
\[\mathcal{B}_q^{\star}=\mathcal{B}_q(qc,a/q,e;\,b,d),\]
\[\mathbb{B}_q^{\star}=\mathbb{B}_q(qc,a/q,e;\,b,d)\]
with $\mathcal{B}_q$ and $\mathbb{B}_q$ being defined in Theorem
\ref{thm-b}. Solving the equation
$\mathcal{B}_q-\mathcal{B}_q^{\star}Y_q=0$ associated with the
variable $Y_q$, we found from equation \eqref{three-hhh} the
following relation.

\begin{thm}[Three-term contiguous relation of $_3\phi_2$-series]\label{thm-aa}
  \bnm
_3\phi_2\ffnk{cccc}{q;\frac{bd}{ace}}{a,c,e}{b,d}
=Y_q\times{_3\phi_2}\ffnk{cccc}{q;\frac{bd}{ace}}{a/q,qc,e}{b,d}+
Z_q\times{_3\phi_2}\ffnk{cccc}{q;\frac{bd}{qace}}{a,qc,qe}{qb,d},
 \enm
where the coefficients $Y_q$ and $Z_q$ are defined by
 \bnm
Y_q= \frac{(a-d)(qce-d)}{(ae-d)(qc-d)},\quad
Z_q=\frac{(1-e)(a-qc)(qace-bd)d}{(1-b)(ae-d)(qc-d)qace}.
 \enm
\end{thm}

Employing the substitutions $a\to q^a$, $b\to q^b$, $c\to q^c$,
$d\to q^d$, $e\to q^e$ for Theorem \ref{thm-aa} and then letting
$q\to1$, we recover the following relation.

\begin{corl}[\citu{wang}{Theorem 42}]
 \bnm
_3F_2\ffnk{cccc}{1}{a,c,e}{b,d}
=Y\,{_3F_2}\ffnk{cccc}{1}{a-1,c+1,e}{b,d}
+Z\,{_3F_2}\ffnk{cccc}{1}{a,c+1,e+1}{b+1,d},
 \enm
where the coefficients $Y$ and $Z$ are given by
 \bnm
Y=\frac{(a-d)(1+c+e-d)}{(a+e-d)(1+c-d)},\quad
Z=\frac{e(1+c-a)(b+d-a-c-e-1)}{b(1+c-d)(a+e-d)}.
 \enm
\end{corl}

Specifying the parameter $b\to a$ in Theorem \ref{thm-aa} and using
$q$-Gauss summation formula \eqref{q-gauss},
 we obtain the following relation with one free parameter less under replacements
 $a\to qa$, $c\to c/q$.

\begin{prop} [Two-term contiguous relation of $_3\phi_2$-series]\label{prop-p}
 \bnm
_3\phi_2\ffnk{cccc}{q;\frac{qd}{ce}}{a,c,e}{qa,d}
&=&{_3\phi_2}\ffnk{cccc}{q;\frac{d}{ce}}{qa,c,qe}{q^2a,d}
\frac{(qa-c)(1-e)d}{(1-qa)(d-qa)ce}\\
&+&\ffnk{cccc}{q}{d/c,d/e}{d,d/ce}_{\infty}\frac{(qae-d)}{(qa-d)e}.
 \enm
\end{prop}

Performing the substitutions $a\to q^a$, $c\to q^c$, $d\to q^d$,
$e\to q^e$ for Proposition \ref{prop-p} and then letting $q\to1$, we
recover the following relation.

\begin{corl}[\citu{wang}{Proposition 43}]
 \bnm
_3F_2\ffnk{cccc}{1}{a,c,e}{a+1,d}&=&{_3F_2}\ffnk{cccc}{1}{a+1,c,e+1}{a+2,d}\frac{e(c-a-1)}{(1+a)(1+a-d)}\\
&+&\Gamma\fnk{cccc}{d,d-c-e}{d-c,d-e}\frac{(1+a+e-d)}{(1+a-d)}.
 \enm
\end{corl}

Taking $d=qae$ in Proposition \ref{prop-p}, we get the following
relation.
\begin{prop} [Two-term contiguous relation of $_3\phi_2$-series]\label{prop-q}
 \bnm
\quad_3\phi_2\ffnk{cccc}{q;\frac{q^2a}{c}}{a,c,e}{qa,qae}
={_3\phi_2}\ffnk{cccc}{q;\frac{qa}{c}}{qa,c,qe}{q^2a,qae}
\frac{(qa-c)}{(qa-1)c}.
 \enm
\end{prop}

Employing the substitutions $a\to q^a$, $c\to q^c$, $e\to q^e$ for
 Proposition \ref{prop-q} and then letting
$q\to1$, we recover the following relation.

\begin{corl}[\citu{wang}{Corollary 44}]
 \bnm
_3F_2\ffnk{cccc}{1}{a,c,e}{a+1,a+e+1}={_3F_2}\ffnk{cccc}{1}{a+1,c,e+1}{a+2,a+e+1}\frac{(1+a-c)}{(1+a)}.
 \enm
\end{corl}

\subsection{\textbf{B}\&\textbf{B}}$ $\\\\
Let Eq$^{\star}$\eqref{eq-b} stand for Eq\eqref{eq-b} under the
parameter replacements
\[b\to d/q, \quad d\to qb.\]
Then for an arbitrary variable $Y_q$, the difference
$\text{Eq}\eqref{eq-b}-Y_q\times\text{Eq}^{\star}\eqref{eq-b}$
results in the relation:
 \bmn
&&_3\phi_2\ffnk{cccc}{q;\frac{bd}{ace}}{a,c,e}{b,d}-Y_q\times
{_3\phi_2}\ffnk{cccc}{q;\frac{bd}{ace}}{a,c,e}{qb,d/q}
  \nnm\\\nnm
&&=\:\big(\mathcal{B}_q-\mathcal{B}_q^{\star}Y_q\big)\times
{_3\phi_2}\ffnk{cccc}{q;\frac{bd}{ace}}{a/q,c,e}{b,d/q}
\\ \label{three-iii}
&&+\:\:\big(\mathbb{B}_q-\mathbb{B}_q^{\star}Y_q\big)\times
{_3\phi_2}\ffnk{cccc}{q;\frac{bd}{qace}}{a,qc,qe}{qb,d},
 \emn
where the following notations have been used for coefficients
\[\mathcal{B}_q^{\star}=\mathcal{B}_q(a,c,e;\,d/q,qb),\]
\[\mathbb{B}_q^{\star}=\mathbb{B}_q(a,c,e;\,d/q,qb)\]
with $\mathcal{B}_q$ and $\mathbb{B}_q$ being defined in Theorem
\ref{thm-b}. Solving the equation
$\mathcal{B}_q-\mathcal{B}_q^{\star}Y_q=0$ associated with the
variable $Y_q$, we derive from equation \eqref{three-iii} the
following relation.

\begin{thm}[Three-term contiguous relation of $_3\phi_2$-series]\label{thm-bb}
  \bnm
_3\phi_2\ffnk{cccc}{q;\frac{bd}{ace}}{a,c,e}{b,d}
=Y_q\times{_3\phi_2}\ffnk{cccc}{q;\frac{bd}{ace}}{a,c,e}{qb,d/q}+
Z_q\times{_3\phi_2}\ffnk{cccc}{q;\frac{bd}{qace}}{a,qc,qe}{qb,d},
 \enm
where the coefficients $Y_q$ and $Z_q$ are defined by
 \bnm
Y_q= \frac{(b-c)(b-e)(q-d)(qce-d)}{(b-1)(b-ce)(qc-d)(qe-d)},\quad
Z_q= \frac{(1-c)(1-e)(d-qb)(qace-bd)}{(b-1)(b-ce)(qc-d)(qe-d)a}.
 \enm
\end{thm}

Performing the substitutions $a\to q^a$, $b\to q^b$, $c\to q^c$,
$d\to q^d$, $e\to q^e$ for Theorem \ref{thm-bb} and then letting
$q\to1$, we recover the following relation.

\begin{corl}[\citu{wang}{Theorem 45}]
 \bnm
_3F_2\ffnk{cccc}{1}{a,c,e}{b,d}
=Y\,{_3F_2}\ffnk{cccc}{1}{a,c,e}{b+1,d-1}
+Z\,{_3F_2}\ffnk{cccc}{1}{a,c+1,e+1}{b+1,d},
 \enm
where the coefficients $Y$ and $Z$ are given by
 \bnm
Y=\frac{(b-c)(b-e)(d-1)(1+c+e-d)}{b(1+c-d)(1+e-d)(c+e-b)},\quad
Z=\frac{ce(1+b-d)(1+a+c+e-b-d)}{b(1+c-d)(1+e-d)(c+e-b)}.
 \enm
\end{corl}

Specifying the parameter $d\to qa$ in Theorem \ref{thm-bb} and using
$q$-Gauss summation formula \eqref{q-gauss},
 we deduce the following relation with one free parameter less.

\begin{prop} [Two-term contiguous relation of $_3\phi_2$-series]\label{prop-r}
 \bnm
_3\phi_2\ffnk{cccc}{q;\frac{qb}{ce}}{a,c,e}{qa,b}
&=&{_3\phi_2}\ffnk{cccc}{q;\frac{b}{ce}}{a,qc,qe}{qa,qb}
\frac{(1-c)(1-e)(a-b)}{(1-b)(a-c)(a-e)}\\
&+&\ffnk{cccc}{q}{b/c,b/e}{b,b/ce}_{\infty}\frac{(1-a)(ce-a)}{(a-c)(a-e)}.
 \enm
\end{prop}

Employing the substitutions $a\to q^a$, $b\to q^b$, $c\to q^c$,
$e\to q^e$ for Proposition \ref{prop-r} and then letting $q\to1$, we
recover the following relation.

\begin{corl}[\citu{wang}{Proposition 46}]
 \bnm
_3F_2\ffnk{cccc}{1}{a,c,e}{a+1,b}&=&{_3F_2}\ffnk{cccc}{1}{a,c+1,e+1}{a+1,b+1}\frac{ce(b-a)}{b(a-c)(a-e)}\\
&+&\Gamma\fnk{cccc}{b,b-c-e}{b-c,b-e}\frac{a(a-c-e)}{(a-c)(a-e)}.
 \enm
\end{corl}

Taking $a=ce$ in Proposition \ref{prop-r}, we attain the following
relation.
\begin{prop} [Two-term contiguous relation of $_3\phi_2$-series]\label{prop-s}
 \bnm
\quad_3\phi_2\ffnk{cccc}{q;\frac{qb}{ce}}{ce,c,e}{qce,b}
={_3\phi_2}\ffnk{cccc}{q;\frac{b}{ce}}{ce,qc,qe}{qce,qb}
\frac{(b-ce)}{(b-1)ce}.
 \enm
\end{prop}

Performing the substitutions $b\to q^b$, $c\to q^c$, $e\to q^e$ for
 Proposition \ref{prop-s} and then
letting $q\to1$, we recover the following relation.

\begin{corl}[\citu{wang}{Corollary 47}]
 \bnm
_3F_2\ffnk{cccc}{1}{c+e,c,e}{c+e+1,b}={_3F_2}\ffnk{cccc}{1}{c+e,c+1,e+1}{c+e+1,b+1}\frac{(b-c-e)}{b}.
 \enm
\end{corl}

\subsection{\textbf{B}\&\textbf{C}}$ $\\\\
Let Eq$^{\star}$\eqref{eq-cc} stand for Eq\eqref{eq-cc} under the
parameter replacements
\[a\to c,\quad c\to a/q^2,\quad e\to e/q,\quad b\to b/q, \quad d\to d/q^2.\]
Then for an arbitrary variable $Y_q$, the difference
$\text{Eq}\eqref{eq-bb}-Y_q\times\text{Eq}^{\star}\eqref{eq-cc}$
leads us to the relation:
 \bmn
&&_3\phi_2\ffnk{cccc}{q;\frac{bd}{ace}}{a,c,e}{b,d}-Y_q\times
{_3\phi_2}\ffnk{cccc}{q;\frac{bd}{ace}}{a/q^2,c,e/q}{b/q,d/q^2}
  \nnm\\\label{three-jj}
&&=\big(\mathcal{B}_q-\mathcal{C}_q^{\star}Y_q\big)\times
{_4\phi_3}\ffnk{cccc}{q;\frac{bd}{qace}}
{a/q,c,e,q\Big(1-\frac{\mathcal{B}_q-\mathcal{C}_q^{\star}Y_q}{\mathfrak{B}_q-\mathfrak{C}_q^{\star}Y_q}\Big)}
{b,d/q,\Big(1-\frac{\mathcal{B}_q-\mathcal{C}_q^{\star}Y_q}{\mathfrak{B}_q-\mathfrak{C}_q^{\star}Y_q}\Big)},
 \emn
where the following notations have been used for coefficients
\[\mathcal{C}_q^{\star}=\mathcal{C}_q(c,a/q^2,e/q;\,b/q,d/q^2),\]
\[\mathfrak{C}_q^{\star}=\mathfrak{C}_q(c,a/q^2,e/q;\,b/q,d/q^2)\]
with $\mathcal{C}_q$ and $\mathfrak{C}_q$ being defined in Theorem
\ref{thm-c}. Solving the equation
$1-\frac{\mathcal{B}_q-\mathcal{C}_q^{\star}Y_q}{\mathfrak{B}_q-\mathfrak{C}_q^{\star}Y_q}
=d/q^2$ associated with the variable $Y_q$, we achieve from equation
\eqref{three-jj} the following relation.

\begin{thm}[Three-term contiguous relation of $_3\phi_2$-series]\label{thm-cc}
 \bnm
_3\phi_2\ffnk{cccc}{q;\frac{bd}{ace}}{a,c,e}{b,d}
=Y_q\times{_3\phi_2}\ffnk{cccc}{q;\frac{bd}{ace}}{a/q^2,c,e/q}{b/q,d/q^2}+
Z_q\times{_3\phi_2}\ffnk{cccc}{q;\frac{bd}{qace}}{a/q,c,e}{b,d/q^2},
 \enm
where the coefficients $Y_q$ and $Z_q$ are defined by
 \bnm
Y_q&=&\frac{(q-b)(q-d)(q^2-d)(bd^2+q^2abce+q^3ace-q^3bce-qabd-qacde)ace}{(q-a)(qc-d)(q^2c-d)(qe-d)(ae-d)qb^2},\\
Z_q&=&\frac{(q-d)(q^2-d)(bd-qace)(abde+q^2bcd+q^3ace-q^3bce-qabd-qacde)}{(q-a)(qc-d)(q^2c-d)(qe-d)(ae-d)qb^2}.
 \enm
\end{thm}

Employing the substitutions $a\to q^a$, $b\to q^b$, $c\to q^c$,
$d\to q^d$, $e\to q^e$ for Theorem \ref{thm-cc} and then letting
$q\to1$, we recover the following relation.

\begin{corl}[\citu{wang}{Theorem 48}]
 \bnm
_3F_2\ffnk{cccc}{1}{a,c,e}{b,d}
=Y\,{_3F_2}\ffnk{cccc}{1}{a-2,c,e-1}{b-1,d-2}
+Z\,{_3F_2}\ffnk{cccc}{1}{a-1,c,e}{b,d-2},
 \enm
where the coefficients $Y$ and $Z$ are given by
 \bnm
Y&=&\frac{(b-1)(d-1)(d-2)(2b+2d+cd+de+2ad-ac-ae-bd-3a-c-e-d^2-1)}{(a-1)(1+c-d)(2+c-d)(1+e-d)(d-a-e)},\\
Z&=&\frac{(d-1)(d-2)(b+d-a-c-e-1)(1+a+c+ae+bd-ad-ce-2b-e)}{(a-1)(1+c-d)(2+c-d)(1+e-d)(d-a-e)}.
 \enm
\end{corl}

Specifying the parameter $d\to qa$ in Theorem \ref{thm-cc} and using
$q$-Gauss summation formula \eqref{q-gauss},
 we establish the following relation with
one free parameter less.

\begin{prop}[Two-term contiguous relation of $_3\phi_2$-series]\label{prop-t}
 \bnm
_3\phi_2\ffnk{cccc}{q;\frac{qb}{ce}}{a,c,e}{qa,b}
&=&{_3\phi_2}\ffnk{cccc}{q;\frac{qb}{ce}}{a/q^2,c,e/q}{a/q,b/q}
\frac{(a-1)(a-q)(a-b)(q-b)c^2e^2}{(a-c)(a-qc)(a-e)(q-e)b^2}\\
&&\xxqdn\xxqdn\xxqdn+\:\ffnk{cccc}{q}{b/c,b/e}{b,b/ce}_{\infty}\frac{(a-1)(ce-b)
(a^2be+q^2abc+q^2ace-q^2bce-qa^2ce-qa^2b)}{(a-c)(a-qc)(a-e)(q-e)b^2}.
 \enm
\end{prop}

Performing the substitutions $a\to q^a$, $b\to q^b$, $c\to q^c$,
$e\to q^e$ for Proposition \ref{prop-t} and then letting $q\to1$, we
recover the following relation.

\begin{corl}[\citu{wang}{Proposition 49}]
 \bnm
_3F_2\ffnk{cccc}{1}{a,c,e}{a+1,b}&=&{_3F_2}\ffnk{cccc}{1}{a-2,c,e-1}{a-1,b-1}
 \frac{a(a-1)(a-b)(b-1)}{(a-c)(a-c-1)(a-e)(e-1)}\\&&\xxqdn
 \qqdn+\:\:\Gamma\fnk{cccc}{b,b-c-e}{b-c,b-e}
 \frac{a(c+e-b)(a^2+ce+e+b-c-ae-ab-1)}{(a-c)(a-c-1)(a-e)(e-1)}.
 \enm
\end{corl}

Taking $b=\frac{qace(a-q)}{a^2e+q^2ac-q^2ce-qa^2}$ in Proposition
\ref{prop-t}, we found the following relation under the replacements
$a\to q^2 a$, $e \to qe$.

\begin{prop}[Two-term contiguous relation of $_3\phi_2$-series]\label{prop-u}
 \bnm
&&{_3\phi_2}\ffnk{cccc}{q;\frac{q^2a(qa-1)}{q^2a^2e+qac-ce-q^2a^2}}{a,c,e}{qa,\frac{qace(qa-1)}{q^2a^2e+qac-ce-q^2a^2}}\\
&=&{_3\phi_2}\ffnk{cccc}{q;\frac{q^2a(qa-1)}{q^2a^2e+qac-ce-q^2a^2}}{q^2a,c,qe}{q^3a,\frac{q^2ace(qa-1)}{q^2a^2e+qac-ce-q^2a^2}}\\
&\times&\frac{(q^2a-c)(qa-e)(qa-1)}{(q^2a-1)(q^2a^2+q^2a^2ce+ce-q^2a^2e-qac-qace)}.
 \enm
\end{prop}

Employing the substitutions $a\to q^a$, $c\to q^c$, $e\to q^e$ for
 Proposition \ref{prop-u} and then letting
$q\to1$, we recover the following relation.

\begin{corl}[\citu{wang}{Corollary 50}]
 \bnm
{_3F_2}\ffnk{cccc}{1}{a,c,e}{a+1,1+a-e+\frac{ce}{a+1}}
&=&{_3F_2}\ffnk{cccc}{1}{a+2,c,e+1}{a+3,2+a-e+\frac{ce}{a+1}}\\
 &\times&\frac{(1+a)(2+a-c)(1+a-e)}{(2+a)(1+2a+a^2+ce-ae-e)}.
 \enm
\end{corl}

\subsection{\textbf{B}\&\textbf{D}}$ $\\\\
Let Eq$^{\star}$\eqref{eq-dd} stand for Eq\eqref{eq-dd} under the
parameter replacements
\[a\to a/q,\quad c\to qc,\quad e\to qe,\quad b\to qb.\]
Then for an arbitrary variable $Y_q$, the difference
$\text{Eq}\eqref{eq-bb}-Y_q\times\text{Eq}^{\star}\eqref{eq-dd}$
results in the relation:
 \bmn
&&{_3\phi_2}\ffnk{cccc}{q;\frac{bd}{ace}}{a,c,e}{b,d}-Y_q\times
{_3\phi_2}\ffnk{cccc}{q;\frac{bd}{ace}}{a/q,qc,qe}{qb,d}
  \nnm\\\label{three-kk}
&&=\big(\mathcal{B}_q-\mathcal{D}_q^{\star}Y_q\big)\times
{_4\phi_3}\ffnk{cccc}{q;\frac{bd}{qace}}
{a/q,c,e,q\Big(1-\frac{\mathcal{B}_q-\mathcal{D}_q^{\star}Y_q}{\mathfrak{B}_q-\mathfrak{D}_q^{\star}Y_q}\Big)}
{b,d/q,\Big(1-\frac{\mathcal{B}_q-\mathcal{D}_q^{\star}Y_q}{\mathfrak{B}_q-\mathfrak{D}_q^{\star}Y_q}\Big)},
 \emn
where the following notations have been used for coefficients
\[\mathcal{D}_q^{\star}=\mathcal{D}_q(a/q,qc,qe;\,qb,d),\]
\[\mathfrak{D}_q^{\star}=\mathfrak{D}_q(a/q,qc,qe;\,qb,d)\]
with $\mathcal{D}_q$ and $\mathfrak{D}_q$ being defined in Theorem
\ref{thm-d}. Solving the equation
$1-\frac{\mathcal{B}_q-\mathcal{D}_q^{\star}Y_q}{\mathfrak{B}_q-\mathfrak{D}_q^{\star}Y_q}
=d/q^2$ associated with the variable $Y_q$, we obtain from equation
\eqref{three-kk} the following relation.

\begin{thm}[Three-term contiguous relation of $_3\phi_2$-series]\label{thm-dd}
 \bnm
_3\phi_2\ffnk{cccc}{q;\frac{bd}{ace}}{a,c,e}{b,d}
=Y_q\times{_3\phi_2}\ffnk{cccc}{q;\frac{bd}{ace}}{a/q,qc,qe}{qb,d}+
Z_q\times{_3\phi_2}\ffnk{cccc}{q;\frac{bd}{qace}}{a/q,c,e}{b,d/q^2},
 \enm
where the coefficients $Y_q$ and $Z_q$ are defined by
 \bnm
Y_q&=&\frac{(a-qb)(a-d)(1-c)(1-e)(q^3ace+q^2abce+bd^2-qabd-q^3bce-qacde)qd}
{(a-q)(1-b)(qc-d)(qe-d)(q^3ace+q^2bcde+bd^2-q^2bcd-q^2bde-qacde)a},\\
Z_q&=&\frac{(q-d)(q^2-d)(qace-bd)(q^2ce+qcde+ad-qcd-qde-qace)}
{(q-a)(qc-d)(qe-d)(q^3ace+q^2bcde+bd^2-q^2bcd-q^2bde-qacde)}.
 \enm
\end{thm}

Performing the substitutions $a\to q^a$, $b\to q^b$, $c\to q^c$,
$d\to q^d$, $e\to q^e$ for Theorem \ref{thm-dd} and then letting
$q\to1$, we recover the following relation.

\begin{corl}[\citu{wang}{Theorem 51}]
 \bnm
_3F_2\ffnk{cccc}{1}{a,c,e}{b,d}
=Y\,{_3F_2}\ffnk{cccc}{1}{a-1,c+1,e+1}{b+1,d}
+Z\,{_3F_2}\ffnk{cccc}{1}{a-1,c,e}{b,d-2},
 \enm
where the coefficients $Y$ and $Z$ are given by
 \bnm
Y&=&\frac{ce(a-b-1)(a-d)\big\{(a-1)(1+c+e-d)+(d-2)(b+d-a-c-e-1)\big\}}
{b(1-a)(1+c-d)(1+e-d)\big\{ce+(d-2)(b+d-a-c-e-1)\big\}},\\
Z&=&\frac{(d-1)(d-2)(1+a+c+e-b-d)(1+c+e+ce+ad-ac-ae-a-d)}
{(1-a)(1+c-d)(1+e-d)\big\{ce+(d-2)(b+d-a-c-e-1)\big\}}.
 \enm
\end{corl}

Specifying the parameter $d\to qa$ in Theorem \ref{thm-dd} and using
$q$-Gauss summation formula \eqref{q-gauss},
 we get the following relation with
one free parameter less.

\begin{prop}[Two-term contiguous relation of $_3\phi_2$-series]\label{prop-v}
 \bnm
_3\phi_2\ffnk{cccc}{q;\frac{qb}{ce}}{a,c,e}{qa,b}
&=&{_3\phi_2}\ffnk{cccc}{q;\frac{qb}{ce}}{a/q,qc,qe}{qa,qb}
\tfrac{(a-b)(a-qb)(1-q)(1-c)(1-e)ce}{(a-c)(a-e)(1-b)(ace+qbc+qbe-qce-qbce-ab)}\\
&+&\ffnk{cccc}{q}{b/c,b/e}{b,qb/ce}_{\infty}
 \frac{ce(1-a)(ace+qac+qae-qce-qace-a^2)}
{(a-c)(a-e)(ace+qbc+qbe-qce-qbce-ab)}.
 \enm
\end{prop}

Employing the substitutions $a\to q^a$, $b\to q^b$, $c\to q^c$,
$e\to q^e$ for Proposition \ref{prop-v} and then letting $q\to1$, we
recover the following relation.

\begin{corl}[\citu{wang}{Proposition 52}]
 \bnm
_3F_2\ffnk{cccc}{1}{a,c,e}{a+1,b}&=&{_3F_2}\ffnk{cccc}{1}{a-1,c+1,e+1}{a+1,b+1}
 \frac{ce(a-b)(1+b-a)}{b(a-c)(a-e)\big\{ce+(a-1)(b-c-e)\big\}}\\
 &+&\Gamma\fnk{cccc}{b,b-c-e+1}{b-c,b-e}
 \frac{a(a^2+ce+c+e-a-ac-ae)}{(a-c)(a-e)\big\{ce+(a-1)(b-c-e)\big\}}.
 \enm
\end{corl}

\subsection{\textbf{B}\&\textbf{D}}$ $\\\\
Let Eq$^{\star}$\eqref{eq-d} stand for Eq\eqref{eq-d} under the
parameter replacements
\[a\to a/q,\quad c\to qc,\quad e\to qe,\quad b\to qb.\]
Then for an arbitrary variable $Y_q$, the difference
$\text{Eq}\eqref{eq-b}-Y_q\times\text{Eq}^{\star}\eqref{eq-d}$ leads
us to the relation:
 \bmn
&&_3\phi_2\ffnk{cccc}{q;\frac{bd}{ace}}{a,c,e}{b,d}-Y_q\times
{_3\phi_2}\ffnk{cccc}{q;\frac{bd}{ace}}{a/q,qc,qe}{qb,d}
  \nnm\\\nnm
&&=\:\big(\mathcal{B}_q-\mathcal{D}_q^{\star}Y_q\big)\times
{_3\phi_2}\ffnk{cccc}{q;\frac{bd}{ace}}{a/q,c,e}{b,d/q}
\\ \label{three-lll}
&&+\:\:\big(\mathbb{B}_q-\mathbb{D}_q^{\star}Y_q\big)\times
{_3\phi_2}\ffnk{cccc}{q;\frac{bd}{qace}}{a,qc,qe}{qb,d},
 \emn
where the following notations have been used for coefficients
\[\mathcal{D}_q^{\star}=\mathcal{D}_q(a/q,qc,qe;\,qb,d),\]
\[\mathbb{D}_q^{\star}=\mathbb{D}_q(a/q,qc,qe;\,qb,d)\]
with $\mathcal{D}_q$ and $\mathbb{D}_q$ being defined in Theorem
\ref{thm-d}. Solving the equation
$\mathbb{B}_q-\mathbb{D}_q^{\star}Y_q=0$ associated with the
variable $Y_q$, we derive from equation \eqref{three-lll} the
following relation.

\begin{thm}[Three-term contiguous relation of $_3\phi_2$-series]\label{thm-ee}
  \bnm
_3\phi_2\ffnk{cccc}{q;\frac{bd}{ace}}{a,c,e}{b,d}
=Y_q\times{_3\phi_2}\ffnk{cccc}{q;\frac{bd}{ace}}{a/q,qc,qe}{qb,d}+
Z_q\times{_3\phi_2}\ffnk{cccc}{q;\frac{bd}{ace}}{a/q,c,e}{b,d/q},
 \enm
where the coefficients $Y_q$ and $Z_q$ are defined by
 \bnm
\quad Y_q=
\frac{(a-qb)(a-d)(1-c)(1-e)qd}{(a-q)(1-b)(qc-d)(qe-d)a},\quad Z_q=
\frac{(q-d)(qace+qcd+qde-qcde-q^2ce-ad)}{(a-q)(qc-d)(qe-d)}.
 \enm
\end{thm}

Performing the substitutions $a\to q^a$, $b\to q^b$, $c\to q^c$,
$d\to q^d$, $e\to q^e$ for Theorem \ref{thm-ee} and then letting
$q\to1$, we recover the following relation.

\begin{corl}[\citu{wang}{Theorem 53}]
 \bnm
_3F_2\ffnk{cccc}{1}{a,c,e}{b,d}
=Y\,{_3F_2}\ffnk{cccc}{1}{a-1,c+1,e+1}{b+1,d}
+Z\,{_3F_2}\ffnk{cccc}{1}{a-1,c,e}{b,d-1},
 \enm
where the coefficients $Y$ and $Z$ are given by
 \bnm
Y=\frac{ce(a-d)(a-b-1)}{b(1-a)(1+c-d)(1+e-d)},\quad
Z=\frac{(d-1)\big\{(a-1)(1+c+e-d)-ce\big\}}{(1-a)(1+c-d)(1+e-d)}.
 \enm
\end{corl}

Specifying the parameter $b\to a$ in Theorem \ref{thm-ee} and using
$q$-Gauss summation formula \eqref{q-gauss},
 we deduce the following relation with one free parameter less under the replacements $a\to qa$, $d\to qd$.

\begin{prop} [Two-term contiguous relation of $_3\phi_2$-series]\label{prop-w}
 \bnm
_3\phi_2\ffnk{cccc}{q;\frac{qd}{ce}}{a,c,e}{qa,d}
&=&{_3\phi_2}\ffnk{cccc}{q;\frac{qd}{ce}}{a,qc,qe}{q^2a,qd}
\frac{(a-d)(1-q)(1-c)(1-e)d}{(1-qa)(d-1)(ace+cd+de-cde-ce-ad)}\\
&+&\ffnk{cccc}{q}{d/c,d/e}{d,qd/ce}_{\infty}\frac{(a-1)ce}{(ace+cd+de-cde-ce-ad)}.
 \enm
\end{prop}

Employing the substitutions $a\to q^a$, $c\to q^c$, $d\to q^d$,
$e\to q^e$ for Proposition \ref{prop-w} and then letting $q\to1$, we
recover the following relation.

\begin{corl}[\citu{wang}{Proposition 54}]
 \bnm
_3F_2\ffnk{cccc}{1}{a,c,e}{a+1,d}&=&{_3F_2}\ffnk{cccc}{1}{a,c+1,e+1}{a+2,d+1}\frac{ce(d-a)}{d(1+a)(ad+ce-ac-ae)}\\
&+&\Gamma\fnk{cccc}{d,d-c-e+1}{d-c,d-e}\frac{a}{(ad+ce-ac-ae)}.
 \enm
\end{corl}

\subsection{\textbf{B}\&\textbf{D}}$ $\\\\
Let Eq$^{\star}$\eqref{eq-d} stand for Eq\eqref{eq-d} under the
parameter replacements
\[a\to c,\quad c\to a,\quad e\to qe,\quad b\to qb.\]
Then for an arbitrary variable $Y_q$, the difference
$\text{Eq}\eqref{eq-b}-Y_q\times\text{Eq}^{\star}\eqref{eq-d}$
results in the relation:
 \bmn
&&_3\phi_2\ffnk{cccc}{q;\frac{bd}{ace}}{a,c,e}{b,d}-Y_q\times
{_3\phi_2}\ffnk{cccc}{q;\frac{bd}{ace}}{a,c,qe}{qb,d}
  \nnm\\\nnm
&&=\:\big(\mathcal{B}_q-\mathcal{D}_q^{\star}Y_q\big)\times
{_3\phi_2}\ffnk{cccc}{q;\frac{bd}{ace}}{a/q,c,e}{b,d/q}
\\ \label{three-mmm}
&&+\:\:\big(\mathbb{B}_q-\mathbb{D}_q^{\star}Y_q\big)\times
{_3\phi_2}\ffnk{cccc}{q;\frac{bd}{qace}}{a,qc,qe}{qb,d},
 \emn
where the following notations have been used for coefficients
\[\mathcal{D}_q^{\star}=\mathcal{D}_q(c,a,qe;\,qb,d),\]
\[\mathbb{D}_q^{\star}=\mathbb{D}_q(c,a,qe;\,qb,d)\]
with $\mathcal{D}_q$ and $\mathbb{D}_q$ being defined in Theorem
\ref{thm-d}. Solving the equation
$\mathcal{B}_q-\mathcal{D}_q^{\star}Y_q=0$ associated with the
variable $Y_q$, we attain from equation \eqref{three-mmm} the
following relation.

\begin{thm}[Three-term contiguous relation of $_3\phi_2$-series]\label{thm-ff}
  \bnm
_3\phi_2\ffnk{cccc}{q;\frac{bd}{ace}}{a,c,e}{b,d}
=Y_q\times{_3\phi_2}\ffnk{cccc}{q;\frac{bd}{ace}}{a,c,qe}{qb,d}+
Z_q\times{_3\phi_2}\ffnk{cccc}{q;\frac{bd}{qace}}{a,qc,qe}{qb,d},
 \enm
where the coefficients $Y_q$ and $Z_q$ are defined by
 \bnm
\quad Y_q= \frac{(b-c)(qce-d)}{(b-1)(qe-d)c},\quad Z_q=
\frac{(1-c)(qace-bd)}{(1-b)(qe-d)ac}.
 \enm
\end{thm}

Performing the substitutions $a\to q^a$, $b\to q^b$, $c\to q^c$,
$d\to q^d$, $e\to q^e$ for Theorem \ref{thm-ff} and then letting
$q\to1$, we recover the following relation.

\begin{corl}[\citu{wang}{Theorem 55}]
 \bnm
3F_2\ffnk{cccc}{1}{a,c,e}{b,d}
=Y\,{_3F_2}\ffnk{cccc}{1}{a,c,e+1}{b+1,d}
+Z\,{_3F_2}\ffnk{cccc}{1}{a,c+1,e+1}{b+1,d},
 \enm
where the coefficients $Y$ and $Z$ are given by
 \bnm
Y=\frac{(b-c)(1+c+e-d)}{b(1+e-d)},\quad
Z=\frac{c(1+a+c+e-b-d)}{b(1+e-d)}.
 \enm
\end{corl}

Specifying the parameter $b\to a$ in Theorem \ref{thm-ff} and using
$q$-Gauss summation formula \eqref{q-gauss},
 we achieve the following relation with one free parameter less under
 the replacement  $e\to e/q$.

\begin{prop}[Two-term contiguous relation of $_3\phi_2$-series]\label{prop-x}
 \bnm
_3\phi_2\ffnk{cccc}{q;\frac{qd}{ce}}{a,c,e}{qa,d}
={_3\phi_2}\ffnk{cccc}{q;\frac{d}{ce}}{a,qc,e}{qa,d}
\frac{(1-c)}{(a-c)}
+\ffnk{cccc}{q}{d/c,d/e}{d,d/ce}_{\infty}\frac{(a-1)}{(a-c)}.
 \enm
\end{prop}

Employing the substitutions $a\to q^a$, $c\to q^c$, $d\to q^d$,
$e\to q^e$ for Proposition \ref{prop-x} and then letting $q\to1$, we
recover the following relation.

\begin{corl}[\citu{wang}{Proposition 56}]
 \bnm
_3F_2\ffnk{cccc}{1}{a,c,e}{a+1,d}={_3F_2}\ffnk{cccc}{1}{a,c+1,e}{a+1,d}\frac{c}{(c-a)}
+\Gamma\fnk{cccc}{d,d-c-e}{d-c,d-e}\frac{a}{(a-c)}.
 \enm
\end{corl}

\subsection{\textbf{C}\&\textbf{C}}$ $\\\\
Let Eq$^{\star}$\eqref{eq-c} stand for Eq\eqref{eq-c} under the
parameter replacements
\[a\to qc,\quad c\to a/q.\]
Then for an arbitrary variable $Y_q$, the difference
$\text{Eq}\eqref{eq-c}-Y_q\times\text{Eq}^{\star}\eqref{eq-c}$ leads
us to the relation:
 \bmn
&&_3\phi_2\ffnk{cccc}{q;\frac{bd}{ace}}{a,c,e}{b,d}-Y_q\times
{_3\phi_2}\ffnk{cccc}{q;\frac{bd}{ace}}{a/q,qc,e}{b,d}
  \nnm\\\nnm
&&=\:\big(\mathcal{C}_q-\mathcal{C}_q^{\star}Y_q\big)\times
{_3\phi_2}\ffnk{cccc}{q;\frac{bd}{ace}}{a,qc,qe}{qb,qd}
\\ \label{three-nnn}
&&+\:\:\big(\mathbb{C}_q-\mathbb{C}_q^{\star}Y_q\big)\times
{_3\phi_2}\ffnk{cccc}{q;\frac{bd}{qace}}{qa,q^2c,q^2e}{q^2b,q^2d},
 \emn
where the following notations have been used for coefficients
\[\mathcal{C}_q^{\star}=\mathcal{C}_q(qc,a/q,e;\,b,d),\]
\[\mathbb{C}_q^{\star}=\mathbb{C}_q(qc,a/q,e;\,q,d)\]
with $\mathcal{C}_q$ and $\mathbb{C}_q$ being defined in Theorem
\ref{thm-c}. Solving the equation
$\mathcal{C}_q-\mathcal{C}_q^{\star}Y_q=0$ associated with the
variable $Y_q$, we establish from equation \eqref{three-nnn} the
following relation.

\begin{thm}[Three-term contiguous relation of $_3\phi_2$-series]\label{thm-gg}
  \bnm
\qquad_3\phi_2\ffnk{cccc}{q;\frac{bd}{ace}}{a,c,e}{b,d}
=Y_q\times{_3\phi_2}\ffnk{cccc}{q;\frac{bd}{ace}}{a/q,qc,e}{b,d}+
Z_q\times{_3\phi_2}\ffnk{cccc}{q;\frac{bd}{qace}}{qa,q^2c,q^2e}{q^2b,q^2d},
 \enm
where the coefficients $Y_q$ and $Z_q$ are defined by
 \bnm
&&Y_q=\frac{(qabce+qacde+bd-abd-qace-qbcde)}{(qabce+qacde+bd-abde-qace-qbcd)},\\
&& Z_q=\frac{(1-a)(a-qc)(1-qc)(1-qe)(1-e)(qace-bd)(bd)^2}
{(1-b)(1-d)(1-qb)(1-qd)(qabce+qacde+bd-abde-qace-qbcd)(qace)^2}.
 \enm
\end{thm}

Performing the substitutions $a\to q^a$, $b\to q^b$, $c\to q^c$,
$d\to q^d$, $e\to q^e$ for Theorem \ref{thm-gg} and then letting
$q\to1$, we recover the following relation.

\begin{corl}[\citu{wang}{Theorem 57}]
 \bnm
_3F_2\ffnk{cccc}{1}{a,c,e}{b,d}
=Y\,{_3F_2}\ffnk{cccc}{1}{a-1,c+1,e}{b,d}
+Z\,{_3F_2}\ffnk{cccc}{1}{a+1,c+2,e+2}{b+2,d+2},
 \enm
where the coefficients $Y$ and $Z$ are given by
 \bnm
Y=\frac{bd-a(1+c+e)}{bd-(1+c)(a+e)},\quad
Z=\frac{ae(1+c)(1+e)(1+c-a)(b+d-a-c-e-1)}{bd(1+b)(1+d)(a+e+ac+ce-bd)}.
 \enm
\end{corl}

Specifying the parameter $d\to qc$ in Theorem \ref{thm-gg} and using
$q$-Gauss summation formula \eqref{q-gauss},
 we found the following relation with one free parameter less.

\begin{prop} [Two-term contiguous relation of $_3\phi_2$-series]\label{prop-y}
 \bnm
_3\phi_2\ffnk{cccc}{q;\frac{qb}{ae}}{a,c,e}{qc,b}
&=&{_3\phi_2}\ffnk{cccc}{q;\frac{b}{ae}}{qa,q^2c,q^2e}{q^3c,q^2b}
\frac{(1-a)(qc-a)(1-e)(1-qe)b^2}{(1-b)(1-qb)(1-qc)(1-q^2c)(ae)^2}\\
&+&\ffnk{cccc}{q}{b/a,b/e}{b,b/ae}_{\infty}\frac{ae(1-qc)-b(1+ae-a-qce)}{e(a-b)(1-qc)}.
 \enm
\end{prop}

Employing the substitutions $a\to q^a$, $b\to q^b$, $c\to q^c$,
$e\to q^e$ for Proposition \ref{prop-y} and then letting $q\to1$, we
recover the following relation.

\begin{corl}[\citu{wang}{Proposition 58}]
 \bnm
_3F_2\ffnk{cccc}{1}{a,c,e}{c+1,b}&=&{_3F_2}\ffnk{cccc}{1}{a+1,c+2,e+2}{c+3,b+2}\frac{ae(1+e)(a-c-1)}{b(1+b)(1+c)(2+c)}
\\&+&\Gamma\fnk{cccc}{b,b-a-e}{b-a,b-e}\frac{(a-b)(1+c)+ae}{(a-b)(1+c)}.
 \enm
\end{corl}

Taking $b=\frac{ae(1-qc)}{1+ae-a-qce}$ in Proposition \ref{prop-y},
we obtain the following relation.

\begin{prop}[Two-term contiguous relation of $_3\phi_2$-series]\label{prop-z}
 \bnm
&&{_3\phi_2}\ffnk{cccc}{q;\frac{q(1-qc)}{1+ae-a-qce}}{a,c,e}{qc,\frac{ae(1-qc)}{1+ae-a-qce}}\\
&=&{_3\phi_2}\ffnk{cccc}{q;\frac{1-qc}{1+ae-a-qce}}{qa,q^2c,q^2e}{q^3c,\frac{q^2ae(1-qc)}{1+ae-a-qce}}\\
&\times&\frac{(1-e)(1-qe)(1-qc)(qc-a)}{(1-q^2c)(1-qce)(1+ae+q^2ace-qae-qce-a)}.
 \enm
\end{prop}

Performing the substitutions $a\to q^a$, $c\to q^c$, $e\to q^e$ for
Proposition \ref{prop-z} and then letting $q\to1$, we recover the
following relation.

\begin{corl}[\citu{wang}{Proposition 59}]
 \bnm
{_3F_2}\ffnk{cccc}{1}{a,c,e}{c+1,a+\frac{ce}{c+1}}
&=&{_3F_2}\ffnk{cccc}{1}{a+1,c+2,e+2}{c+3,2+a+\frac{ce}{c+1}}\\
 &\times&\frac{(a-c-1)(1+c)(1+e)e}{(2+c)(1+c+e)(1+c+a+ac+ae)}.
 \enm
\end{corl}
\subsection{\textbf{D}\&\textbf{D}}$ $\\\\
Let Eq$^{\star}$\eqref{eq-d} stand for Eq\eqref{eq-d} under the
parameter replacements
\[a\to c/q,\quad c\to qa.\]
Then for an arbitrary variable $Y_q$, the difference
$\text{Eq}\eqref{eq-d}-Y_q\times\text{Eq}^{\star}\eqref{eq-d}$
results in the relation:
 \bmn
&&_3\phi_2\ffnk{cccc}{q;\frac{bd}{ace}}{a,c,e}{b,d}-Y_q\times
{_3\phi_2}\ffnk{cccc}{q;\frac{bd}{ace}}{qa,c/q,e}{b,d}
  \nnm\\\nnm
&&=\:\big(\mathcal{D}_q-\mathcal{D}_q^{\star}Y_q\big)\times
{_3\phi_2}\ffnk{cccc}{q;\frac{bd}{ace}}{a,c/q,e/q}{b/q,d/q}
\\ \label{three-ooo}
&&+\:\:\,\big(\mathbb{D}_q-\mathbb{D}_q^{\star}Y_q\big)\times
{_3\phi_2}\ffnk{cccc}{q;\frac{bd}{qace}}{qa,c,e}{b,d},
 \emn
where the following notations have been used for coefficients
\[\mathcal{D}_q^{\star}=\mathcal{D}_q(c/q,qa,e;\,q,d),\]
\[\mathbb{D}_q^{\star}=\mathbb{D}_q(c/q,qa,e;\,b,d)\]
with $\mathcal{D}_q$ and $\mathbb{D}_q$ being defined in Theorem
\ref{thm-d}. Solving the equation
$\mathbb{D}_q-\mathbb{D}_q^{\star}Y_q=0$ associated with the
variable $Y_q$, we get from equation \eqref{three-ooo} the following
relation.

\begin{thm}[Three-term contiguous relation of $_3\phi_2$-series]\label{thm-hh}
  \bnm
_3\phi_2\ffnk{cccc}{q;\frac{bd}{ace}}{a,c,e}{b,d}
=Y_q\times{_3\phi_2}\ffnk{cccc}{q;\frac{bd}{ace}}{qa,c/q,e}{b,d}+
Z_q\times{_3\phi_2}\ffnk{cccc}{q;\frac{bd}{ace}}{a,c/q,e/q}{b/q,d/q},
 \enm
where the coefficients $Y_q$ and $Z_q$ are defined by
 \bnm
\quad Y_q= \frac{(a-1)(c-b)(c-d)aq^2}{(c-q)(qa-b)(qa-d)c},\quad Z_q=
\frac{(q-b)(q-d)(qa-c)a}{(q-c)(qa-b)(qa-d)}.
 \enm
\end{thm}

Employing the substitutions $a\to q^a$, $b\to q^b$, $c\to q^c$,
$d\to q^d$, $e\to q^e$ for Theorem \ref{thm-hh} and then letting
$q\to1$, we recover the following relation.

\begin{corl}[\citu{wang}{Theorem 60}]
 \bnm
_3F_2\ffnk{cccc}{1}{a,c,e}{b,d}
=Y\,{_3F_2}\ffnk{cccc}{1}{a+1,c-1,e}{b,d}
+Z\,{_3F_2}\ffnk{cccc}{1}{a,c-1,e-1}{b-1,d-1},
 \enm
where the coefficients $Y$ and $Z$ are given by
 \bnm
Y=\frac{a(b-c)(c-d)}{(1-c)(1+a-b)(1+a-d)},\quad
Z=\frac{(1-b)(1-d)(1+a-c)}{(1-c)(1+a-b)(1+a-d)}.
 \enm
\end{corl}

Instead, solving the equation
$\mathcal{D}_q-\mathcal{D}_q^{\star}Y_q=0$ associated with the
variable $Y_q$, we derive from equation \eqref{three-ooo} the
following relation.

\begin{thm}[Three-term contiguous relation of $_3\phi_2$-series]\label{thm-ii}
  \bnm
_3\phi_2\ffnk{cccc}{q;\frac{bd}{ace}}{a,c,e}{b,d}
=Y_q\times{_3\phi_2}\ffnk{cccc}{q;\frac{bd}{ace}}{qa,c/q,e}{b,d}+
Z_q\times{_3\phi_2}\ffnk{cccc}{q;\frac{bd}{qace}}{qa,c,e}{b,d},
 \enm
where the coefficients $Y_q$ and $Z_q$ are defined by
 \bnm
\quad Y_q= \frac{(c-b)(c-d)aq}{(qa-b)(qa-d)c},\quad Z_q=
\frac{(qa-c)(qace-bd)}{(qa-b)(qa-d)ce}.
 \enm
\end{thm}

Performing the substitutions $a\to q^a$, $b\to q^b$, $c\to q^c$,
$d\to q^d$, $e\to q^e$ for Theorem \ref{thm-ii} and then letting
$q\to1$, we recover the following relation.

\begin{corl}[\citu{wang}{Theorem 61}]
 \bnm
3F_2\ffnk{cccc}{1}{a,c,e}{b,d}
=Y\,{_3F_2}\ffnk{cccc}{1}{a+1,c-1,e}{b,d}
+Z\,{_3F_2}\ffnk{cccc}{1}{a+1,c,e}{b,d},
 \enm
where the coefficients $Y$ and $Z$ are given by
 \bnm
Y=\frac{(c-b)(c-d)}{(1+a-b)(1+a-d)},\quad
Z=\frac{(1+a-c)(1+a+c+e-b-d)}{(1+a-b)(1+a-d)}.
 \enm
\end{corl}



\end{document}